\documentclass{autart}
\usepackage{amsmath,amssymb}
\usepackage{array}
\usepackage{mdwmath}
\usepackage{mdwtab}
\usepackage[mathscr]{eucal}
\usepackage{graphicx}
\usepackage{times}
\usepackage{pifont}
\usepackage{srcltx}
\usepackage{subfigure}
\usepackage{color}
\usepackage{epsfig}
\usepackage{enumerate}

\newcommand{\ov}[1]{\overline{#1}}

\newcommand{\bldr}{\boldsymbol{r}}

\newcommand{\col}{\,\mbox{\bf col}\,}
\newcommand{\bt}{\boldsymbol{\mathscr{T}}}

\newcommand{\sgn}{\text{\rm sgn}}
\newcommand{\url}[1]{\texttt{#1}}
\renewcommand{\Box}{\bullet}
\newcommand{\dsafe}{d_{\text{safe}}}
\newcommand{\dtrig}{d_{\text{trig}}}
\newcommand{\drange}{d_{\text{range}}}
\newcommand{\dist}[1]{\text{\bf dist}_D[#1]}
\newcommand{\appr}{^{\text{\tiny $\approx$}}\! 0}
\newcommand{\fr}{\text{\large $\mathfrak{r}$}}
\newcommand{\bro}{\boldsymbol{\rho}}

\newcommand{\sphan}[2]{\sphericalangle_{#1} #2}
\newcommand{\tang}[1]{\sphericalangle\text{\bf TANG}\left[#1\right]}

\newcommand{\fd}{\mathfrak{d}}

\newcommand{\brro}{\boldsymbol{\varrho}}

\renewcommand{\thefootnote}{\fnsymbol{footnote}}
\newtheorem{remark}{Remark}
\newtheorem{definition}{Definition}
\newtheorem{assumption}{Assumption}
\newtheorem{corollary}{Corollary}
\newtheorem{observation}{Observation}

\newcommand{\epf}{$\quad \bullet$}

\begin{document}
\begin{frontmatter}
\title{Proofs of the Technical Results Justifying a Biologically Inspired Algorithm for Reactive Navigation of Nonholonomic Robots in Maze-Like Environments}
\author[auth1]{Alexey S. Matveev}\ead{almat1712@yahoo.com},
\author[auth2]{Michael C. Hoy}\ead{mch.hoy@gmail.com},
\author[auth2]{Andrey V. Savkin}\ead{a.savkin@unsw.edu.au}

\address[auth1]{Department of Mathematics and Mechanics, Saint Petersburg
University, Universitetskii 28, Petrodvoretz, St.Petersburg, 198504,
Russia}

\address[auth2]{School of Electrical Engineering and
Telecommunications, The University of New South Wales, Sydney 2052,
Australia}
\end{frontmatter}
\section{Introduction}
Inspired by behaviors of animals, which are believed to use
simple, local motion control rules that result in remarkable
and complex intelligent behaviors \cite{Thomp66,CaJh99,Faj01}, we examine the navigation strategy that is aimed at reaching a steady target in a steady arbitrarily shaped maze-like environment and is composed of the following reflex-like rules:
\begin{enumerate}[{\bf s.1)}]
\item At considerable distances from the obstacle,
\begin{enumerate}
\item
\label{rule.a}
turn towards the target as quickly as possible;
\item move directly to the target when headed to it;
\end{enumerate}
\item At a short distance from the obstacle,
\begin{enumerate}
\setcounter{enumii}{2}
\item Follow (a,b) when leaving from the obstacle;
\item When approaching it, quickly avert the collision threat by sharply turning.
\end{enumerate}
\end{enumerate}
 Studies of target pursuit in animals, ranging from dragonflies to fish and dogs to humans, have suggested that they often use the pure pursuit guidance s.1)
 to catch not only a steady but also a moving target. The idea of local obstacle avoidance strategy s.2) is also inspired by biological examples such as a cockroach encountering a wall \cite{CaJh99}.
\par
The rules s.1), s.2) demand only minor perceptual capacity. Access even to the distance to the obstacle is not needed: it suffices to determine whether it is short or not, and be aware of the sign of its time derivative. As for the target, the vehicle has to access its relative bearing angle. Moreover, it suffices that it is able only to recognize which quadrant of its relative Cartesian frame hosts the target line-of-sight.
\par
\par
To address the issue of nonholonomic constraints, control saturation, and under-actuation, we consider a vehicle of the Dubins car type. It
is capable of moving with a constant speed along planar paths of
upper limited curvature without reversing the direction and is controlled by the upper limited angular velocity. As a result, it is unable to slow down, stop, or make an abrupt turn.
\par
By reliance on the bearing-only data about the target, the proposed approach is similar to the Pledge algorithm \cite{AbSe80} and Angulus algorithm \cite{LuTi91}. Unlike ours, the both assume access to the absolute direction (e.g., by a compass), and the latter employs not one but two angles in the convergence criterion. The major distinction is that they assume the vehicle to be able to trace the paths of unlimited curvature, in particular, broken curves and to move exactly along the obstacle boundary. These assumptions are violated in the context of this paper, which entails deficiency in the available proofs of the convergence of these algorithms.
\par
The extended introduction and discussion of the proposed control law are given in the paper submitted by the authors to the IFAC journal Automatica.
This text basically contains the proofs of the technical facts underlying justification of the convergence at performance of the proposed algorithm in that paper, which were not included into it due to the length limitations. To make the current text logically consistent, were reproduce the problem statement and notations.
\section{Problem Setup and the Navigation Strategy}
\label{S2}
We consider a planar under-actuated nonholonomic vehicle of the Dubins car type. It  travels with a constant speed $v$ without reversing direction
and is controlled by the angular velocity $u$ limited by a given
constant $\ov{u}$. There also is a steady point target $\bt$ and a single steady
obstacle $D \not\ni \bt$ in the plane, which is an arbitrarily shaped compact domain whose boundary $\partial D$ is Jordan piece-wise analytical curve without inner corners. Modulo smoothened approximation of such  corners, this assumption is typically satisfied by all obstacles encountered in robotics, including continuous mazes. The objective is to drive the vehicle to the target
with constantly respecting a given safety margin $d(t) \geq d_{\text{safe}} >0$. Here $d(t)$ is the distance to the obstacle
\begin{equation}
\label{distance} d(t) := \text{\bf dist}_{D}[\bldr(t)], \qquad
\dist{\bldr}:= \min_{\bldr_\ast \in D} \|\bldr_\ast - \bldr\|,
\end{equation}
$\|\cdot\|$ is the Euclidian norm, and $\bldr(t)$ is the vehicle position.
\par
This position is given by the abscissa $x$ and ordinate $y$ of the vehicle in the world frame, whereas its orientation is described by the angle $\theta$ from the abscissa axis to the robot centerline. The kinematics of the considered vehicles are classically described by the following equations:
\begin{equation}
\label{1}
\begin{array}{l}
\dot{x} = v \cos \theta,
\\
\dot{y} = v \sin \theta,
\end{array}, \quad\dot{\theta} = u \in [-\overline{u},
\overline{u}], \quad
\begin{array}{l}
\bldr(0) = \bldr_0 \not\in D
\\
\theta(0) = \theta_0
\end{array} .
\end{equation}
Thus the minimal turning radius of the vehicle is equal to
\begin{equation}
\label{Rmin} R= v/\overline{u}.
\end{equation}
The vehicle has access to
the current distance $d(t)$ to $D$ and the sign $\sgn \dot{d}(t)$ of its time-rate
$\dot{d}(t)$,
 which are accessible only within the given sensor range: $d \leq d_{\text{range}}$, where $
d_{\text{range}} > d_{\text{safe}}$. The vehicle also has
access to the angle $\beta$ from its forward centerline ray to the
target.
\par
To specify the control strategy s.1), s.2), we introduce the threshold $\dtrig < d_{\text{range}}$ separating the 'short' and 'long' distances to the obstacle.
Mathematically, the examined strategy is given by the following concise formula:
\begin{gather}
\label{c.a}
u = \overline{u} \times
\begin{cases}
\sgn \beta \hfill \mid & \!\!\! \text{if}\; d >\dtrig \; (\text{mode}\,\mathfrak{A})
\\
\left\{
\begin{array}{l l}
\sgn \beta & \text{if} \; \dot{d} > 0
\\
- \sigma & \text{if} \; \dot{d} \leq 0
\end{array} \right| & \!\!\! \text{if}\; d \leq \dtrig \; (\text{mode}\,\mathfrak{B})
\end{cases}.
\end{gather}
Here $\sigma=\pm$ is a constant controller parameter, which gives the turn direction, and $\dot{d}\geq 0$ and $\dot{d}<0$ are equivalent to the vehicle orientation outwards and towards $D$.
The switch $\mathfrak{A} \mapsto \mathfrak{B}$ occurs when
$d$ reduces to $\dtrig$; the converse switch holds when
$d$ increases to $\dtrig$.
When mode $\mathfrak{B}$ is activated, $\dot{d} \leq 0$; if $\dot{d}=0$, the 'turn' submode $u:=-\sigma \ov{u}$ is set up.
Since the control law \eqref{c.a} is discontinuous, the solution of
the closed-loop system is meant in the Fillipov's sense
\cite{UT92}.
\begin{remark}
\rm
\label{rem1}
In \eqref{c.a}, $\beta$  accounts for not only the  heading but also the sum of full turns performed by the target bearing. 
\end{remark}
\par
In the {\it basic} version of the algorithm, the parameter $\sigma$ is fixed. To find a target hidden deeply inside the maze, a modified version can be employed:
whenever $\mathfrak{A} \mapsto \mathfrak{B}$, the parameter $\sigma$ is updated. The updated value is picked randomly and independently of the previous choices from
$\{+, -\}$, with the value $+$ being drawn with a fixed probability
$p \in (0,1)$. This version is called the {\it randomized} control law.
\par
To state the assumptions, we introduce the Frenet frame
$T(\bldr_\ast), N(\bldr_\ast)$ of $\partial D$ at the point $\bldr_\ast \in \partial D$ ($T$ is the positively oriented unit tangent vector, $N$ is the unit normal vector directed inwards $D$, the boundary is oriented so that when traveling on $\partial D$ one has $D$ to the left),
$\varkappa(\bldr_\ast )$ is the
signed curvature ($\varkappa(\bldr_\ast)
< 0$ on concavities) and $R_{\varkappa}(\bldr_\ast) :=
|\varkappa(\bldr_\ast)|^{-1}$. Due to the absence of inner corners, any point $\bldr \not\in D$ at a sufficiently small distance $\text{\bf dist}_D[\bldr] < d_\star$ from $D$ does not belong to the focal locus of $\partial D$ and $\text{\bf dist}_D[\bldr]$ is attained at only one point \cite{Krey91}.
The {\em regular margin} $d_\star(D)>0$ of $D$ is the supremum of such $d_\star$'s.
So $d_\star(D) = \infty$ for convex domains; for non-convex $D$,
\begin{equation}
\label{uniq.rad}
d_\star (D) \leq R_D := \inf_{\bldr \in \partial D: \varkappa(\bldr) <0} R_\kappa(\bldr) .
\end{equation}
(The infimum over the empty set is set to be $+\infty$.)
\begin{assumption}
\label{ass.d}
The vehicle is maneuverable enough: it is capable of full turn without violation of a safety margin $\dsafe > R$ within the regularity margin of the maze $3R < d_\star(D)$, and moreover $4R < R_D$.
\end{assumption}
\begin{assumption}
\label{ass.r}
The sensor range gives enough space to avoid collision with $D$ after its detection: $d_{\text{range}} > 3R$.
\end{assumption}
\par
The parameters $d_{\text{trig}}$ and $\dsafe$ are tuned so that
\begin{equation}
\label{eq.trig} 3R< \dsafe + 2R< \dtrig < d_\star(D) , d_{\text{range}},
R_D-R.
\end{equation}
Such a choice is possible thanks to Assumptions~\ref{ass.d} and \ref{ass.r}.
\section{Main Results}
\label{S5}
\begin{thm}
\label{th.main}
{\bf (i)}
With probability $1$, the randomized control law drives the vehicle at the target $\bt$ for a finite time with always respecting the safety margin (i.e., there exists a time instant $t_\ast$ such that  $\bldr(t_\ast) = \bt$ and
$\dist{\bldr(t)} \geq \dsafe \; \forall t \in [0,t_\ast]$) whenever
both the vehicle initial location $\bldr_0$ and the target are far enough from
  the obstacle and from each other:
\begin{equation}
\label{far.enough} \dist{\bldr_0}> \dtrig + 2R,  \|\bldr_0-\bt\|>2R, \dist{\bt} >
\dtrig .
\end{equation}
\par
{\bf (ii)} The basic control law drives the vehicle at the target for a finite time with always respecting the safety margin whenever \eqref{far.enough} holds and the vehicle initial location and the target lie far enough from the convex hull $\text{\bf co}\;D$ of the maze: $\text{\bf dist}_{\text{\bf co}\;D}[\bt] > \dtrig, \text{\bf dist}_{\text{\bf co}\;D}[\bldr_0] > \dtrig$.
\end{thm}
In \eqref{far.enough}, $\dist{\bldr_0}> \dtrig + 2R$ can be relaxed into $\dist{\bldr_0} >
\dtrig$ if the vehicle is initially directed to the target
$\beta(0)=0$. In view of \eqref{Rmin} and the freedom \eqref{eq.trig} in the choice of $\dsafe, \dtrig$, not only Assumptions~\ref{ass.d}, \ref{ass.r} but also the constraints \eqref{far.enough} disappear (are boiled down into $\dist{\bldr_0}> 0 ,  \|\bldr_0-\bt\|> 0, \dist{\bt} >
0 $) as $v \to 0$. In other words, the algorithm succeeds in any case if the cruise speed $v$ is small enough.
\par
The last assumption $\text{\bf dist}_{\text{\bf co}\;D}[\bt] > \dtrig$ from {\bf (ii)} can be relaxed to cover some scenarios with the target inside the maze. To specify this, we need some notations and definitions.
\par
The $d$-{\it equidistant curve} $C(d)$ of $D$ is the locus of points
$\bldr$ at the distance $\text{\bf dist}_D[\bldr] = d$ from $D$; the $d$-{\it neighborhood} $N(d)$ of $D$ is the area bounded by $C(d)$;
$[\bldr_1,\bldr_2]$ is the straight line segment directed from $\bldr_1$ to $\bldr_2$.
\par
Let $\bldr_\lozenge, \bldr_\ast \in C(\dtrig)$ and $(\bldr_\lozenge, \bldr_\ast)\cap N(d_{\text{trig}}) = \emptyset$.
The points $\bldr_\lozenge, \bldr_\ast$ divide $C(\dtrig)$ into
two arcs. Being concatenated with $[\bldr_\lozenge, \bldr_\ast]$,
each of them gives rise to a Jordan curve encircling a bounded
domain, one of which is the other united with $N(\dtrig)$. The smaller
domain is called the {\it simple cave of $N(\dtrig)$ with endpoints} $\bldr_\lozenge, \bldr_\ast$.
The location $\bldr$ is said to be {\it locked} if it belongs to a simple cave of $N(\dtrig)$ whose endpoints lie on a common ray centered at $\bt$. We remark that if $\text{\bf dist}_{\text{\bf co}\;D}[\bldr] > \dtrig$, the location is unlocked.
\begin{thm}
\label{th.maina}
The basic control law drives the vehicle at the target for a finite time with always respecting the safety margin whenever \eqref{far.enough} holds and both the initial location of the vehicle and the target are unlocked.
\end{thm}
\par
Now we disclose the tactical behavior implied by s.1), s.2) and
show that it includes wall following in a sliding mode.
In doing so, we focus on a particular {\it avoidance maneuver}
 ({\bf AM}), i.e., the motion within uninterrupted mode
 $\mathfrak{B}$.
 \par
Let $\bro(s)$ be the natural parametric representation of $\partial
D$, where $s$ is the curvilinear abscissa. This abscissa is cyclic: $s$ and
$s+L$ encode a common point, where $L$ is the perimeter of $\partial
D$. We notationally identify $s$ and $\bro(s)$. For any $\bldr
\not\in D$ within the regular margin $\dist{\bldr} < d_\star(D)$,
the symbol $s(\bldr)$ stands for the boundary point closest to
$\bldr$, and $s(t):= s[\bldr(t)]$, where $\bldr(t)$ is the vehicle
location at time $t$.
\par
To simplify the matters, we first show that $\partial D$ can be assumed $C^1$-smooth without any loss of generality. Indeed, if
$0 < d < d_\star(D)$, the equidistant curve $C(d)$ is $C^1$-smooth and piece-wise $C^2$-smooth
\cite{Krey91}; its parametric representation, orientation, and
curvature are given by
\begin{equation}
\label{eqdc} s \mapsto \bro(s) - d N(s), \qquad \varkappa_{C(d)}(s)
= \frac{\varkappa(s)}{1+\varkappa(s)d}.
\end{equation}
The second formula holds if $s$ is not a corner point of $\partial
D$; such points contribute circular arcs of the radius $d$ into
$C(d)$. So by picking $\delta>0$ small enough, expanding $D$ to $N(\delta)$, and correction $\fd := \fd - \delta$ of $\fd:=d, \dsafe, \dtrig, \drange$, we keep all assumptions true and do not alter the operation of the closed-loop system.
 Hence $\partial D$ can be assumed $C^1$-smooth.
  \par
  Writing $f(\eta_\ast \pm \appr) >0$ means that there exists small enough $\Delta >0$ such that $f(\eta) >0$ if $0 < \pm (\eta- \eta_\ast )<  \Delta$. The similar notations, e.g., $f(\eta_\ast \pm \appr) \leq 0$, are defined likewise.
  \begin{prop}
\label{lem.cir} Let for the vehicle driven by the control law
\eqref{c.a}, obstacle avoidance be started with zero target bearing
$\beta(t)=0$ at $t=t_\ast$. Then the following claims hold:
\begin{enumerate}[(i)]
\item
\label{init.item} There exists $\tau \geq t_\ast$ such that the
vehicle moves with the maximal steering angle $u \equiv - \sigma \ov{u}
$ and the distance to the obstacle decreases $\dot{d}\leq 0$
until $\tau$,\footnote{This part of AM is called the {\it initial
turn} and abbreviated IT.} and at $t=\tau$, the sliding motion along
the equidistant curve $C\left\{\dist{\bldr(\tau)}\right\}$
\footnote{This is abbreviated SMEC and means following the wall at the fixed distance  $\dist{\bldr(\tau)}$, which is set up at the start of SMEC.} is started with
$\sigma \dot{s}>0$ and $\beta \dot{s}
>0$;
\item \label{rule1} SMEC holds until $\beta$ arrives at
zero at a time when $\varkappa [s(t)+\sigma \appr ] > 0$, which
sooner or later holds and after which a straight move to the
target\footnote{SMT, which is sliding motion over the surface $\beta=0$} is commenced;
\item \label{rule2}
During SMT, the vehicle first does not approach the obstacle
$\dot{d}\geq 0$ and either the triggering threshold $\dtrig$ is
ultimately trespassed and so mode $\mathfrak{B}$ is switched off, or
a situation is encountered where $\dot{d}(t)=0$ and $\varkappa
[s(t)+\sigma \appr ] < 0$. When it is encountered, the vehicle
starts SMEC related to the current distance;
    \item There may be several transitions from SMEC to SMT and
    vice versa, all obeying the rules from (\ref{rule1}), (\ref{rule2});
    \item The number of transitions is finite and finally the vehicle does
    trespass the triggering threshold $\dtrig$, thus terminating
    the considered avoidance maneuver;
 \item \label{dir.item} Except for the initial turn described in (\ref{init.item}), the vehicle maintains a definite direction of bypassing the obstacle: $\dot{s}$ is constantly positive if $\sigma=+$ (counterclockwise bypass) and negative if $\sigma=-$ (clockwise bypass).
 \end{enumerate}
\end{prop}
By \eqref{c.a}, AM is commenced with $\dot{d}(t_\ast)\leq 0$.
The next remark shows that if $\dot{d}(t_\ast)=0$, IT may have the zero
duration.
\begin{remark}
\rm \label{rem.dotzero} If $\dot{d}(t_\ast)=0$, IT has the zero
duration if and only if $\sigma \dot{s}(t_\ast)>0$. Then the
following claims are true:
\begin{enumerate}
\item If $\varkappa[s(t_\ast) + \sigma \cdot \appr] < 0$,
SMEC is immediately started;
\item If $\varkappa[s(t_\ast) + \sigma \cdot \appr] \geq 0$,
the duration of SMEC is zero, and SMT is continued.
\end{enumerate}
\end{remark}
\par
The assumption $\beta(t_\ast)=0$ of Proposition~\ref{lem.cir} holds
for the first AM due to \eqref{far.enough}. Indeed, since
$\dist{\bldr_0}>\dtrig+2R$, the pursuit guidance law turns the
vehicle towards the target earlier than the threshold $\dtrig$ for
activation of AM is encountered. It also holds for all subsequent
AM's since any AM ends in course of SMT by
Proposition~\ref{lem.cir}.
\section{Technical facts underlying the proofs of Proposition~\ref{lem.cir} and
Remark~\ref{rem.dotzero}.} \label{sec.prpr}
\subsection{Geometrical Preliminaries}
\begin{figure}
\subfigure[]{\scalebox{0.25}{\includegraphics{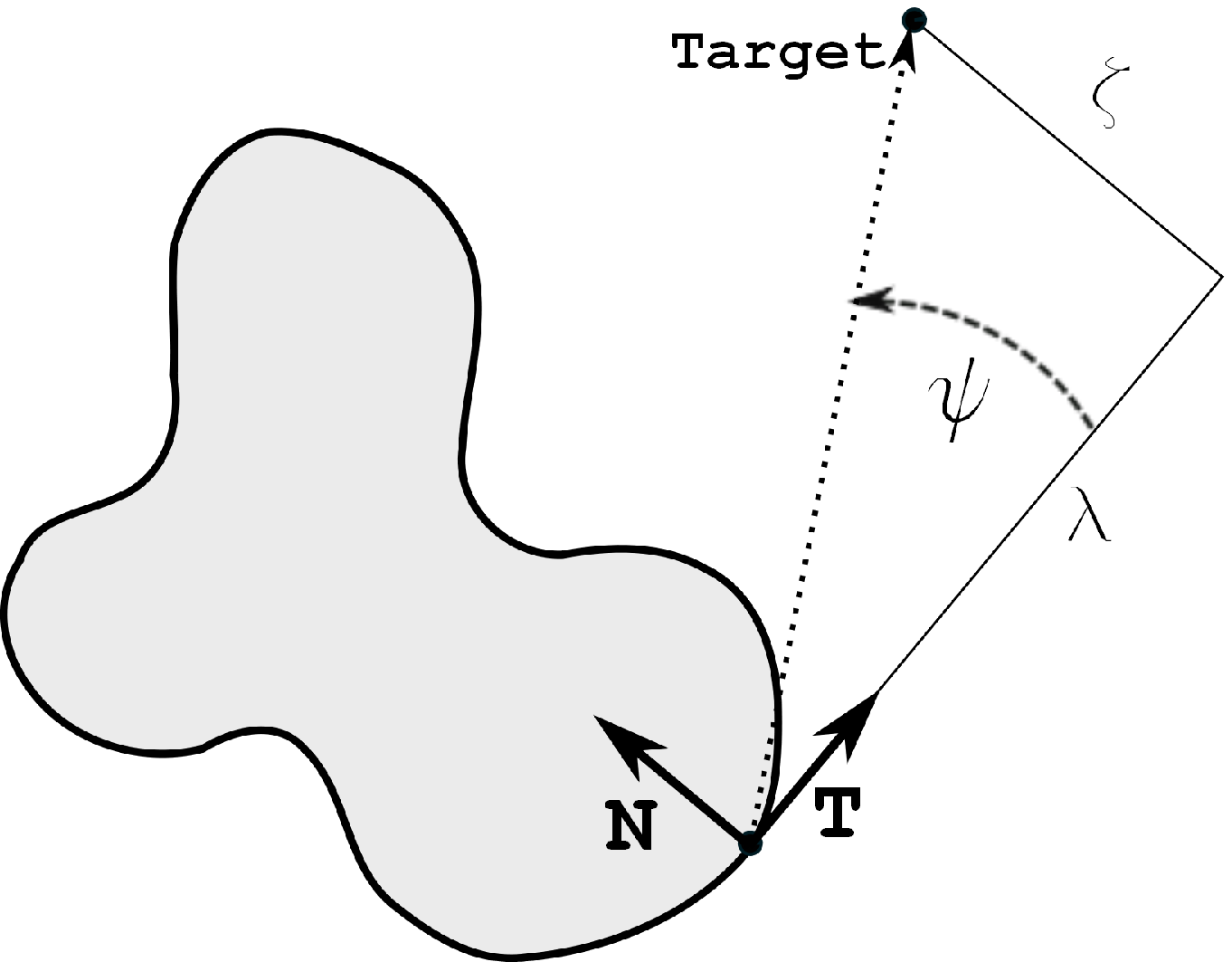}}}
\subfigure[]{\scalebox{0.25}{\includegraphics{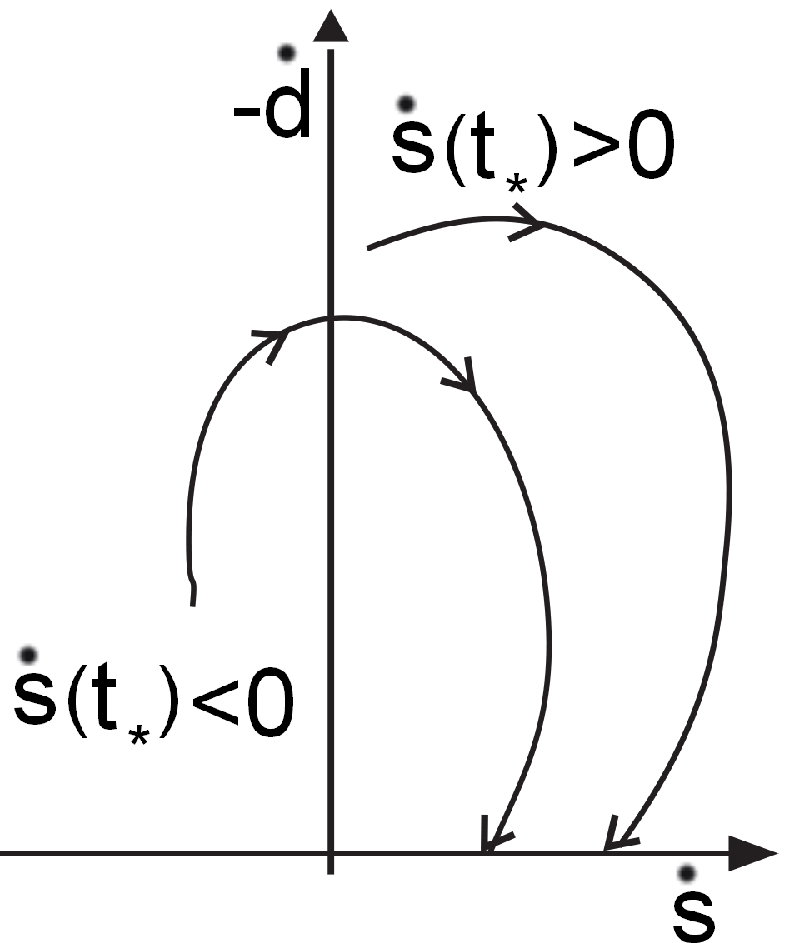}}}
\caption{(a) Definition of $\lambda$ and $\zeta$; (b) Behavior
during IT.} \label{fig.newc}
\end{figure}

We assume that the world frame (WF) is centered at the target $\bt$.
Let $C \not\ni \bt$ be a regular piece-wise smooth directed curve
with natural parametric representation $\brro(s), s \in [s_-,s_+]$.
The turning angle of $C$  around a point $\boldsymbol{p} \not\in C$
is denoted by $\sphericalangle_{\boldsymbol{p}} C$, and $\tang{C}:=
\sphan{0}{T}$, where $T(s), N(s)$ is the Frenet frame of $C$ at
$s$.\footnote{\label{foot.conv} At the corner points, the count of
$\sphan{0}{T}$ progresses abruptly according to the conventional
rules \cite{Krey91}.}
Let $\lambda(s), \zeta(s)$ and $\psi(s)$ stand for the Cartesian
coordinates and polar angle of $- \brro(s)$ in this frame (see
Fig.\ref{fig.newc}(a)), respectively,
and let $^\prime$ denote differentiation with respect to $s$. The
polar angle of $\brro(s)$ in WF and the curvature of $C$ at $s$ are
denoted by $\varphi(s)$ and $\varkappa(s)$, respectively. To
indicate the curve $C$, the symbols $T,N,\lambda,\zeta, \varkappa$,
etc. may be supplied with the lower index $_C$. The directed curve
traced as $s$ runs from $s_1$ to $s_2$ is denoted by
$C_{s_1\xrightarrow{\pm} s_2}$, where the specifier $\pm$ is used
for closed curves. The superscript $^a$ means that the lemma is equipped with the number under which its formulation is given in the basic version of the paper.
\begin{lem}\hspace{-5.0pt}$^a$
The following relations hold whenever $\bt \not\in
C$:
\begin{gather}
\label{geom.eq}
\begin{array}{l}
\lambda^\prime = -1 + \varkappa \zeta \\
 \zeta^\prime = - \varkappa
\lambda
\end{array},
\quad
\begin{array}{r}
\psi^\prime = - \varkappa  + \zeta(\lambda^2+\zeta^2)^{-1}
\\
\varphi^\prime = \zeta(\lambda^2+\zeta^2)^{-1} \end{array},
\\
\label{eq.angle} \fr := \col(\lambda,\zeta) \neq 0, \sphan{0}{\fr} =
\sphan{\bt}{C} - \tang{C}.
\end{gather}
\end{lem}
\pf Differentiation of the equation $\bt = \brro+\lambda T +\zeta N$
and the Frenet-Serret formulas $T^\prime = \kappa N , N^\prime = -
\kappa T$ \cite{Krey91} yield that $ 0 = T + \lambda^\prime T +
\lambda \varkappa N + \zeta^\prime N - \zeta \varkappa T. $ Equating
the cumulative coefficients in this linear combination of $T$ and
$N$ to zero gives the first two equations in \eqref{geom.eq}. By
virtue of them, the third and forth ones follow from \cite{Krey91}
\begin{equation}
\label{angle.rot} \psi^\prime = \frac{\zeta^\prime \lambda -
\lambda^\prime \zeta}{\lambda^2+\zeta^2}, \qquad \varphi^\prime =
\frac{y^\prime x - x^\prime y}{x^2+y^2}.
\end{equation}
The first relation in \eqref{eq.angle} holds since $\bt \not\in C$.
Let $\eta(s):= \sphericalangle \text{\bf TANG} \big[ T_{s_- \to s-0}
\big] + \eta_0$, where $\eta_0$ is the polar angle of $T(s_-)$. The
matrix $\Phi_{\eta(s)}$ of rotation through $\eta(s)$ transforms the
world frame into the Frenet one, and $\brro(s) = h(s) \col[\cos
\varphi(s), \sin \varphi(s)]$. So $\fr(s) = -\Phi_{- \eta(s)}
\brro(s) = h(s) \col\{[\cos [\pi+ \varphi(s) - \eta(s)], \sin
[\pi+\varphi(s) - \eta(s)]\}$. Thus $\pi+ \varphi(s) - \eta(s)$ is
the piece-wise continuous polar angle of $\fr(s)$ that jumps
according to the convention concerned by footnote
$^{\text{\ref{foot.conv}}}$. This trivially implies
\eqref{eq.angle}. \epf
\begin{corollary}
\label{cor.rotat}
 Let $\zeta(s_\ast)=0$ and $\varsigma = \pm$. Then
\begin{equation}
\label{eq.cross}
\begin{array}{r c l}
\varsigma \zeta[s_\ast + \varsigma \appr] \sgn \lambda[s_\ast] <0
&\text{\rm if}& \varkappa[s_\ast+\varsigma \appr] >0
\\
\varsigma \zeta[s_\ast + \varsigma \appr] \sgn \lambda[s_\ast] > 0
&\text{\rm if}& \varkappa[s_\ast+\varsigma \appr] < 0
\end{array}.
\end{equation}
\end{corollary}
\begin{figure}
\subfigure[]{\scalebox{0.24}{\includegraphics{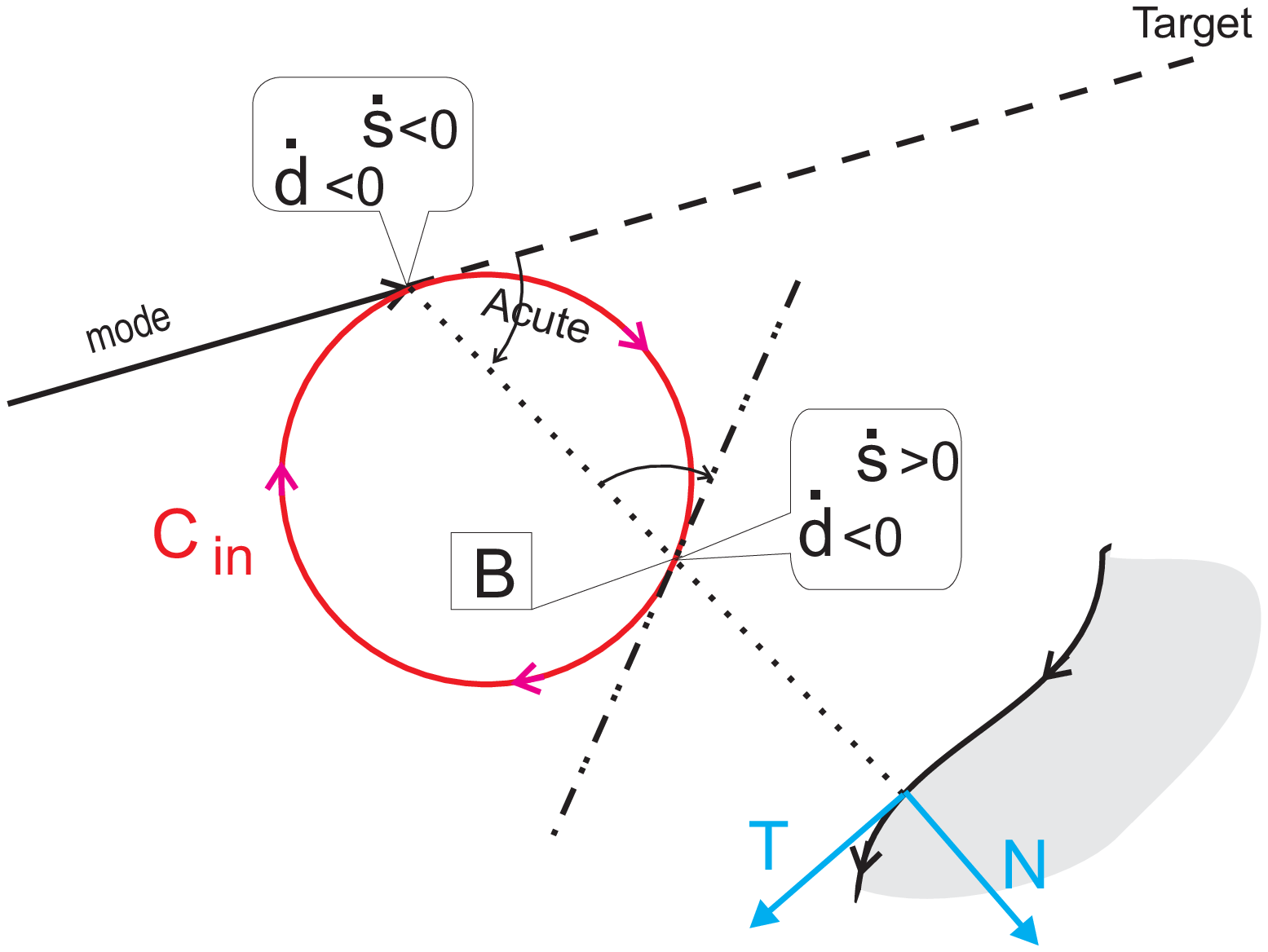}}}
\subfigure[]{\scalebox{0.2}{\includegraphics{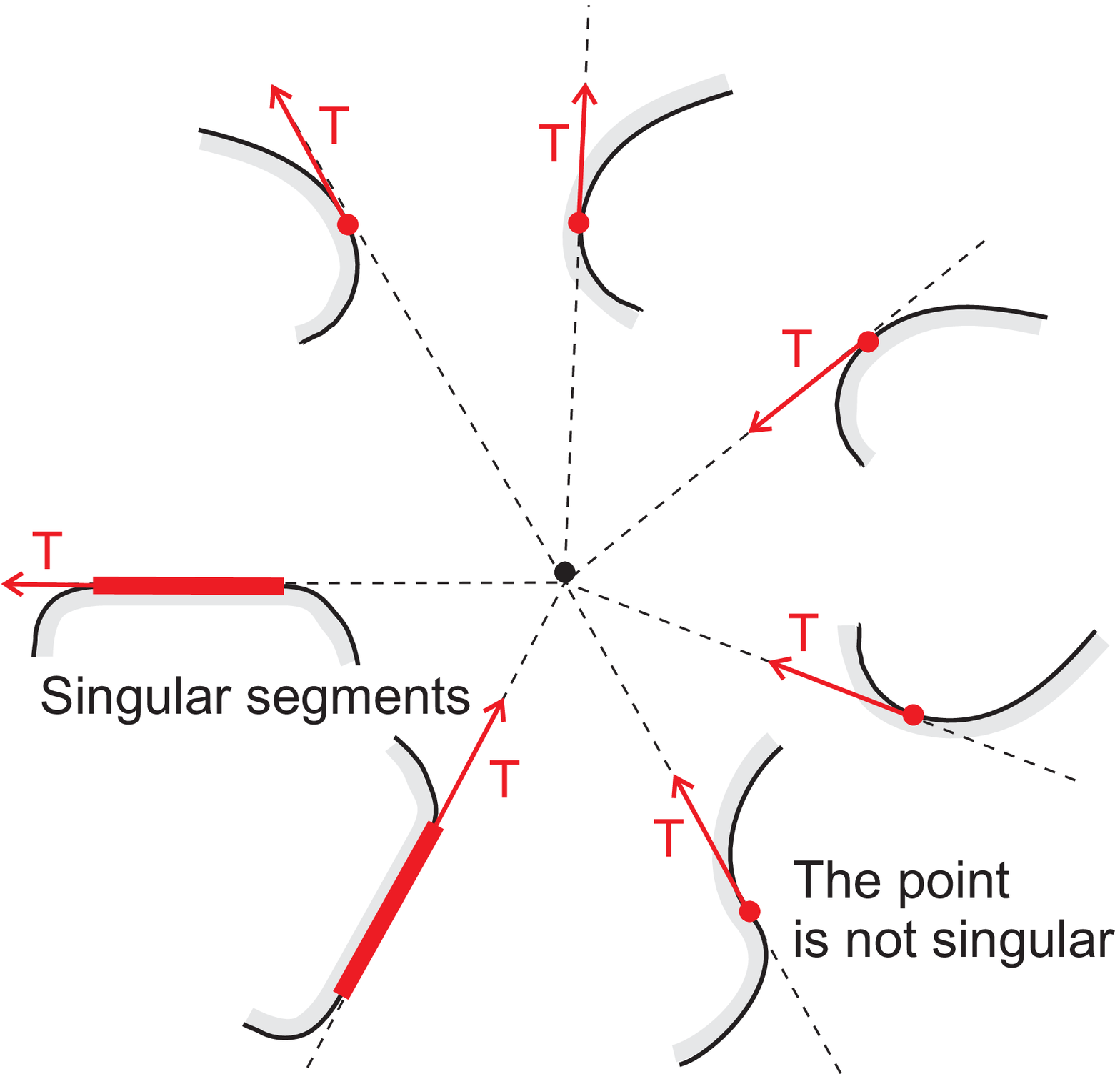}}}
\caption{ (a) Behavior during IT; (b) Singular points.}
\label{fig.cross}
\end{figure}
By \eqref{eq.trig} and the last inequality in \eqref{far.enough},
Lemma~\ref{lem.geom} yields
\begin{equation}
\label{rem.geom} \sphan{0}{\fr_{C(d_\ast)}} = - 2 \pi \qquad
\text{for} \quad d_\ast \in [0,\dtrig].
\end{equation}
\begin{corollary}
\label{cor.finite} There exist $F$ and $d_{\#}> \dtrig$ such that
whenever $|d| \leq d_{\#}$, the set $S(d):= \{s \in
\partial D: \zeta_{\partial D}(s) = d \}$ has no more than $F$ connected components.
\end{corollary}
\pf By the last inequality in \eqref{far.enough}, $\exists d_{\#}:
\dtrig < d_{\#} < \dist{\bt} \leq \sqrt{\zeta(s)^2 + \lambda(s)^2}$.
Then $s \in S(d) \wedge |d|\leq d_{\#} \Rightarrow |\lambda(s)| \geq
\delta := \sqrt{\dist{\bt}^2-d_{\#}^2} >0$. Since the domain $D$ is
compact, $|\lambda^\prime(s)| \leq M < \infty \; \forall s$. So
whenever $s \in S(d)$ and $|d| \leq d_{\#}$, the function
$\lambda(\cdot)$ does not change its sign in the $\delta
M^{-1}$-neighborhood $V(s)$ of $s$.
\par
Since $\partial D$ is piece-wise analytical, each set $\{s: \pm \varkappa(s)>0\}$
and $\{s: \varkappa(s)=0\}$ has finitely many connected components
$\partial^\pm_i$ and $\partial^0_\nu$, respectively. By the
foregoing and \eqref{geom.eq}, any intersection $V(s)\cap
\overline{\partial^\pm_i}, s \in S(d), |d|\leq d_{\#}$ contains only
one point $s$. Hence the entire arc $\partial^\pm_i$ of the length
$\left| \partial^\pm_i \right|$ contains no more than
$\delta^{-1}M\left|
\partial^\pm_i \right|+1$ such points. It remains to note that $S(d)$
covers any $\partial^0_\nu$ such that $\partial^0_\nu \cap S(d) \neq
\emptyset$. \epf
\begin{observation}
\label{obs.exit} SMEC with $\sigma = \pm$ ends when $s \in S_0:= \{s
\in
\partial D: -d_{\#}<\zeta_{\partial D}(s) <0,
\pm \lambda_{\partial D}(s) >0\}$. This set has no more than $F$
connected components, called $\pm${\em arcs}.
\end{observation}
The second claim holds since $\lambda^\prime <0$ on $S_0$ due to
\eqref{eq.trig}, \eqref{geom.eq}.
\subsection{Technical Facts}
\begin{lem}
\label{lem.slide} The following two statements hold:
\\
{\bf (i)}
In the domain $d \leq \dtrig \wedge \dot{d}>0 \vee d > \dtrig$, the
surface $\beta = 0$ is sliding, with the equivalent control {\rm \cite{UT92}} $u\equiv
0$;
\\
{\bf (ii)}
The surface $\dot{d}=0$ is sliding in the domain
\begin{equation}
\label{sub-dom} \dtrig - 2R \leq d < \dtrig, \quad \dot{s}\beta>0,
\quad \sigma \dot{s} >0.
\end{equation}
\end{lem}
\pf {\bf (i)} Let $h$ be the distance from the vehicle to $\bt$. Due
to \eqref{1}, $\dot{h} = - v \cos
\beta,\dot{\beta} = h^{-1} v \sin \beta -u$.
So as the state approaches the surface $\beta =0$,
we have
$\dot{\beta} \overset{\eqref{c.a}}{\to} - \overline{u}\sgn \beta$,
which implies the first claim. 
\par
{\bf (ii)} Let $\alpha$ be the polar angle of the vehicle velocity
in the frame $T_{\partial D}[s(t)], N_{\partial D}[s(t)]$. By
\eqref{uniq.rad}, \eqref{eq.trig}, and \eqref{sub-dom}, $1+
\varkappa[s(t)] d(t) >0$, and as is shown in e.g., \cite{MaHoSa10},
\begin{equation}
\label{eq1} \dot{s} = \frac{v \cos\alpha }{1+ \varkappa(s) d }, \;
\dot{d} =  - v \sin \alpha, \; \dot{\alpha} = - \varkappa(s) \dot{s}
+u .
\end{equation}
As the state approaches a point where $\dot{d}=0$ and
\eqref{sub-dom} holds,
\begin{equation}
\label{bas.eq}
\begin{array}{l}
\sin \alpha \to 0
\\
\cos \alpha \to \sgn \dot{s} \end{array}, \ddot{d} \to - v^2 \left[
\frac{u}{v} \sgn \dot{s} - \frac{\varkappa}{1+\varkappa d}\right].
\end{equation}
If the state remains at a definite side of the surface $\dot{d}=0$,
\eqref{Rmin} and \eqref{c.a} yield that
\begin{gather}
\nonumber \ddot{d} \xrightarrow{\dot{d}>0} \ddot{d}_+:= - v^2 \left[
\frac{1}{R} \sgn (\beta \dot{s}) - \frac{\varkappa}{1+\varkappa
d}\right]
\\ \nonumber
\overset{\text{\eqref{sub-dom}}}{=} - v^2 \left[ \frac{1}{R} -
\frac{\varkappa}{1+\varkappa d}\right], \qquad \ddot{d}
\xrightarrow{\dot{d}< 0} \ddot{d}_- :=
\\
\label{ddots} v^2 \left[  \sigma\frac{1}{R} \sgn \dot{s} +
\frac{\varkappa}{1+\varkappa d}\right]
\overset{\text{\eqref{sub-dom}}}{=} v^2 \left[  \frac{1}{R} +
\frac{\varkappa}{1+\varkappa d}\right].
\end{gather}
The proof is completed by observing that by \eqref{eq.trig},
\eqref{sub-dom},
\begin{gather}
\ddot{d}_+ =  - v^2 \frac{1 + \varkappa d - \varkappa
R}{R(1+\varkappa d)} <0 \text{since}
\begin{array}{l}
1+\varkappa d >0 \; \text{and}
\\
d > \dsafe >R
\end{array}
\nonumber
\\
\label{ddminus} \ddot{d}_- =  v^2 \frac{|\varkappa|
\left[R_\varkappa + (d + R)\sgn\varkappa \right]}{R(1+\varkappa d)}
> 0 .
\end{gather}
\par
The subsequent proofs are focused on $\sigma=+$; the case $\sigma=-$
is considered likewise.
\begin{lem}
If $\dot{d}(t_\ast)<0$, claim (i) in Proposition~{\rm \ref{lem.cir}}
is true.
 \label{lem.i}
\end{lem}
\pf Let $\sigma=+$. Due to \eqref{c.a}, initially $u \equiv -
\ov{u}$. Let $[t_\ast, \tau]$ denote the maximal interval on which
$u \equiv - \ov{u}$. For $t \in (t_\ast,\tau)$, the vehicle moves
clockwise along a circle $C_{\text{in}}$ of the radius $R$ and so by
Remark~\ref{rem1}, $\beta(t)>0$ and
\begin{multline}
\label{far.enough1} d(t) \geq \dist{\bldr(0)} - \underbrace{\|\bldr
- \bldr(0)\|}_{\leq 2R}
\\
\geq \dtrig -2R \overset{\text{\eqref{eq.trig}}}{>} \dsafe >R >0.
\end{multline}
\begin{multline*}
\dot{\alpha} \overset{\text{\eqref{eq1}}}{=} - v \left[  \frac{
\varkappa \cos \alpha}{1+ \varkappa d } +  \frac{\ov{u}}{v} \right]
\\
 \overset{\text{\eqref{Rmin}}}{\leq} - v \left[ \frac{1}{R} -
\frac{ |\varkappa|}{1+ \varkappa d }  \right] = - v \left[
\frac{1}{R} - \frac{1}{R_\varkappa + d \sgn \varkappa}  \right] .
\end{multline*}
While $d \leq \dtrig$ (in particular, while $\dot{d} \leq 0$) the
expression in the last square brackets is positive. This is true by
\eqref{far.enough1} if $\varkappa \geq 0$; otherwise, since
$R_{\varkappa} >R+\dtrig$ by \eqref{eq.trig}. So $\dot{\alpha} \leq
- \delta <0$, i.e., the vector $\col(\cos \alpha, \sin \alpha)$
rotates clockwise. Here the signs of the first and second components
equal those of $\dot{s}$ and $-\dot{d}$, respectively, by
\eqref{eq1} and so $\col(\dot{s},\dot{d})$ evolves as is illustrated
in Fig.~\ref{fig.newc}(b). This and the conditions \eqref{sub-dom}
for the sliding motion complete the proof. \epf
\par
More can be derived from the above proof.
\setcounter{thm}{8}
\begin{lem}\hspace{-5.0pt}$^a$
 Let $s_\ast$ and $s_b$ be the values of the
continuously evolving $s$ at the start and end of IT, respectively.
During IT, $\sigma \dot{s} \geq 0$ if $\sigma \dot{s}(t_\ast) \geq
0$, and $\dot{s}$ ones changes the sign otherwise. In any case, $s$
runs from $s_\ast$ to $s_b$ in the direction $\sigma$ during a last
phase of IT.
\end{lem}
\setcounter{thm}{6}
\pf Let $\sigma = +$.
  The map $\bldr \mapsto (s,d)$ is the
  orientation-reversing immersion on the disc $D_{\text{in}}$
encircled by $C_{\text{in}}$. So it transforms any negatively
oriented circle $C \subset D_{\text{in}}$ concentric with
$C_{\text{in}}$ into a curve $\xi$ with $\tang{\xi} = 2 \pi$. Then
the argument from the concluding part of the proof of
Lemma~\ref{lem.i} shows that as the robot once runs over
$C_{\text{in}}$ in the negative direction, the vector
$\col(\dot{s},\dot{d})$ intersects the half-axes of the frame in the
order associated with counter clockwise rotation, each only once.
This immediately implies the claim given by the first sentence in
the conclusion of the lemma.
  \par
  If $\dot{s}(t_\ast) \geq 0$, this claim yields
that $s_b-s_\ast \geq 0$. Let  $\dot{s}(t_\ast) < 0$. As the robot
once runs over $C_{\text{in}}$ in the negative direction,
$\dot{s}>0$ and $\dot{d} \leq 0$ when it passes the point $B$ from
Fig.~\ref{fig.cross}(a), which corresponds to the second passage of
$s=s_\ast$. Due to the order in which $\col(\dot{s},\dot{d})$
intersects the half-axes, this combination of signs is possible only
before $\dot{d}$ vanishes for the first time, i.e., within IT. Thus
the second occurrence of $s=s_\ast$ holds within IT. The proof is
completed by noting that $\dot{s} >0$ after this by the first claim
of the lemma. \epf
 \par
We proceed to the case where some of the vector fields is tangential
to the discontinuity surface $\dot{d}=0$. Since this may undermine
uniqueness of the solution (its existence is still guaranteed), the
arguments become much more sophisticated. The first lemma
establishes a required technical fact. To state it, we note that
whenever $d:=\dist{\bldr} < R_\star(D)$, the system state
$(x,y,\theta)$ is given by $s,d,\theta$ and along with
$(\dot{d},\dot{s}) \neq (0,0)$, uniquely determines $\beta \in
(-\pi, \pi)$.
\begin{lem}
If $\lambda_{C(d_\dagger)}(s_\ast)\neq 0$ for $d_\dagger \in
[0,\dtrig]$, there exists $\delta>0$ such that whenever $s_\ast \leq
s_0 < s < s_\ast+\delta$ and $|d_\ast - d_\dagger|< \delta$, the
following entailments hold with $\varsigma:= \sgn \dot{s}$:
\begin{gather}
\label{beta.1}
\begin{array}{c}
\dot{s}  \neq 0, \dot{d} \geq 0, d \geq d_\ast,
\zeta_{C(d_\ast)}(s_0) \geq 0
\\
\varkappa(s_\ast+ \varsigma\appr)<0 , \dot{s}
\lambda_{C(d_\dagger)}(s_0)
> 0
\end{array}
\Rightarrow \dot{s} \beta >0;
\\
\begin{array}{c}
\dot{s}  \neq 0, \dot{d} \leq 0, d \leq d_\ast,
\zeta_{C(d_\ast)}(s_0) \leq 0
\\
\varkappa(s_\ast+ \varsigma\appr)\geq 0 , \dot{s}
\lambda_{C(d_\dagger)}(s_0)
> 0
\end{array}
\Rightarrow \dot{s} \beta \leq 0. \label{beta.2}
\end{gather}
In \eqref{beta.2}, $\dot{s} \beta < 0$ if $\zeta_{C(d_\ast)}(s_0) <
0$ or $\varkappa \not\equiv 0$ on $\partial D_{s_0 \to s}$.
\end{lem}
\pf We pick $\delta>0$ so that $\lambda_{C(d_\ast)}(s)$ and
$\varkappa(s)$ do not change the sign as $s$ and $d_\ast$ run over
$(s_\ast, s_\ast+\delta)$ and $(d_\dagger- \delta, d_\dagger+
\delta)$, respectively. By \eqref{eqdc}, the curvature
$\varkappa_{C(d_\ast)}(s)$ does not change its sign either, which
equals $\sgn \varkappa(s_\ast+ \varsigma \appr)$.
\par
If the conditions from \eqref{beta.1} hold and $\varsigma = +$,
application of the second equation from \eqref{geom.eq} to
$C(d_\ast)$ yields that $\zeta_{C(d_\ast)}(s) > 0$. So the target
polar angle in the $s$-related Frenet frame of $C(d_\ast)$ belongs
to $(0,\pi/2)$. Transformation of this frame into that of the
vehicle path consists in a move of the origin in the negative
direction along the $\zeta$-axis (since $d\geq d_\ast$) and a
clockwise rotation of the axes (since $\dot{d}>0, \dot{s}>0$). Since
both operations increase the target bearing angle, $\beta>0$.
Formula \eqref{beta.1} with $\varsigma = -$ and \eqref{beta.2} are
established likewise. \epf
\setcounter{thm}{6}
\begin{lem}\hspace{-5.0pt}$^a$
 Let $\dsafe\leq d_\ast:=d(t_\ast) \leq \dtrig$, $\dot{d}(t_\ast) =0$
 at a time $t_\ast$ within mode $\mathfrak{B}$.
Then for $t>t_\ast, t \approx t_\ast$, the robot
\begin{enumerate}[{\bf i)}]
\item  performs the turn with $u \equiv - \sigma \ov{u}$
if $\sigma \dot{s}(t_\ast) <0$, $d(t_\ast) = \dtrig$, and
$\beta(t_\ast)=0$;
    \item undergoes SMEC if $\sigma \dot{s}(t_\ast) > 0$ and either
   {\bf (1)} $\sigma \beta(t_\ast)>0$ or {\bf (2)} $\beta(t_\ast) =0$ and
    $\varkappa[s(t_\ast) + \sgn \dot{s}(t_\ast) \appr] < 0$;
  \item moves straight to the target if
  $\beta(t_\ast)=0, \sigma \dot{s}(t_\ast) > 0, \varkappa[s(t_\ast) + \sgn \dot{s}(t_\ast) \appr] \geq 0$.
\end{enumerate}
\end{lem}
\setcounter{thm}{8}
\pf Let $\sigma =+$. {\bf i)} As $t \to t_\ast$, \eqref{c.a} and
\eqref{bas.eq} yield that
\begin{equation}
\label{u.minus} \ddot{d}|_{u=-\ov{u}} \to v^2 \left[ - \frac{1}{R} +
\frac{\varkappa}{1+\varkappa d_\ast}\right] = - \frac{1+\varkappa
[d_\ast - R]]}{R(1+\varkappa d_\ast)} <0,
\end{equation}
where $\varkappa:= \varkappa[s(t_\ast)\pm 0]$ and the inequality
holds since $d_\ast \geq \dsafe >R$ due to \eqref{eq.trig}.
\par
Let i) fail to be true and $\varkappa[s(t_\ast) - \appr] < 0$. If
there exists an infinite sequence $\{t_i\}$ such that $t_i>t_\ast,
d(t_i) < \dtrig \; \forall i$ and $t_i \to t_\ast$ as $i \to
\infty$, a proper decrease of every $t_i$ yields in addition that
$\dot{d}(t_i)<0$ since $d(t_\ast)= \dtrig$. However then $\dot{d}(t)
 < 0 $ for $t\geq t_i, t \approx t_\ast$ by \eqref{c.a},
\eqref{u.minus} and thus
 $\dot{d}(t)<0, d(t)< \dtrig$ for $t>t_\ast, t \approx t_\ast$,
 i.e., (i) holds in violation of the initial assumption. It follows
 that $d(t_\ast+\appr) \geq \dtrig$.
\par
Now suppose that there is a sequence $\{t_i\}$ such that
$t_i>t_\ast, d(t_i) = \dtrig \; \forall i$, $t_i \to t_\ast$ as $i
\to \infty$. Then $\dot{d}(t_i)=0$ and so $\beta(t_i)<0$ due to
\eqref{beta.1}. By continuity, $\beta <0$ in a vicinity of the
system state at $t=t_i$. Then any option from \eqref{c.a} yields
$u=-\ov{u}$ and so $u(t) \equiv - \ov{u}\; \forall t \approx t_i$ by
the definition of Filippov's solution. Hence $d(t_i) = \dtrig \wedge
\dot{d}(t_i)=0\overset{\text{\eqref{u.minus}}}{\Rightarrow}
d(t_i+\appr) < \dtrig$, in violation of the foregoing. So $d
> \dtrig$ and $u=\sgn \beta$ for $t>t_\ast, t \approx t_\ast$ by
\eqref{c.a}, and by Lemma~\ref{lem.slide}, SMT is continued. Then
the last relation in \eqref{bas.eq} (with $u:=0$) and
$\varkappa[s(t_\ast) - \appr] < 0$ imply the contradiction
$d(t_\ast+\appr) < \dtrig$ to the foregoing, which proves i).
\par
Let $\varkappa[s(t_\ast) - \appr] \geq 0$. So far as the controller
is first probationally set to the submode related with $\dot{d}<0$,
this submode will be maintained longer by \eqref{u.minus}.
\par
{\bf ii.1)} If $d(t_\ast) < \dtrig$, the claim is true by
Lemma~\ref{lem.slide}. Let $d(t_\ast) = \dtrig$. If there is a
sequence $\{t_i\}$ such that $t_i>t_\ast, d(t_i) < \dtrig \; \forall
i$ and $t_i \to t_\ast$ as $i \to \infty$, a proper decrease of
every $t_i$ yields in addition that $\dot{d}(t_i)<0$. Let $\tau_i$
be the minimal $\tau \in [t_\ast,t_i]$ such that $d(t)< \dtrig$ and
$\dot{d}(t)<0$ for $t \in (\tau,t_i]$. For such $t$, $u \equiv -
\ov{u}$ by \eqref{c.a} and so $\ddot{d}>0$ by \eqref{ddots} and
\eqref{ddminus}. So $\dot{d}(\tau_i) < \dot{d}(t_i) <0, \tau_i >
t_\ast$, and $d(\tau_i) = \dtrig$, otherwise $\tau_i$ is not the
minimal $\tau$. Thus at time $\tau_i$, the assumptions of
Lemma~\ref{lem.i} hold except for $\beta(\tau_i) =0$. In the proof
of this lemma, this relation was used only to justify that $\beta
>0$, which is now true by assumption and the continuity argument. So
by Lemmas~\ref{lem.slide} and \ref{lem.i}, sliding motion along an
equidistant curve $C(d_\dagger)$ with $d_\dagger < \dtrig$ is
commenced at the time $t> \tau_i$ when $\dot{d}(t)=0$ and maintained
while $\beta>0$ and $\dot{s}>0$, in violation of $d(\tau_i) = \dtrig
\; \forall i \wedge \tau_i \xrightarrow{i \to \infty} t_\ast$. This
contradiction proves that $d(t_\ast+\appr) \geq 0$.
\par
Now suppose that there exists a sequence $\{t_i\}$ such that
$t_i>t_\ast, d(t_i) > \dtrig \; \forall i$ and $t_i \to t_\ast$ as
$i \to \infty$. Since $d(t_\ast)=0$, a proper perturbation of every
$t_i$ yields in addition that $\dot{d}(t_i) >0$. Let $\tau_i$ be the
minimal $\tau \in [t_\ast,t_i]$ such that $d(t)> \dtrig$ for $t \in
(\tau,t_i]$. For such $t$, the continuity argument gives $\beta>0$,
\eqref{c.a} yields $u \equiv \ov{u}$ and so $\ddot{d}< 0$ by
\eqref{ddots} and \eqref{ddminus}. Hence $\dot{d}(\tau_i)>0, \tau_i
> t_\ast, d(\tau_i) = \dtrig$ and so $d(\tau_i-\appr)<0$, in
violation of the foregoing. This contradiction proves that
$d(t_\ast+\appr) \equiv 0$ indeed.
\par
 {\bf ii.2)} We first assume that $d_\ast < \dtrig$.
 Due to \eqref{ddots} and \eqref{ddminus}
 \begin{equation}
 \label{ddot.plus}
 \ddot{d}|_{u = -\ov{u}} > 0 \quad\text{and} \quad \ddot{d}|_{u = \ov{u}} < 0 \quad \text{for} \quad t \approx
 t_\ast.
 \end{equation}
So it is easy to see that $\dot{d}(t_\ast+\appr) \geq 0$ and
$d(t_\ast+\appr) \geq d_\ast $. Suppose that $\dot{d}(t_\ast+\appr)
\not\equiv 0$ and so $d(t_\ast+\appr) >d_\ast$. In any
right-vicinity $(t_\ast, t_\ast+\delta)$, there is $\tau$ such that
$\dot{d}(\tau)>0$. For any such $\tau$ that lies sufficiently close
to $t_\ast$, \eqref{beta.1} yields $\beta(\tau)>0$. So $u = \ov{u}$
by \eqref{c.a} and $\ddot{d}(\tau)<0$ by \eqref{ddot.plus}. Hence
the inequality $\dot{d}(t)>0$ is not only maintained but also
enhanced as $t$ decreases from $\tau$ to $t_\ast$, in violation of
the assumption $\dot{d}(t_\ast)=0$ of the lemma. This contradiction
shows that $\dot{d}(t_\ast+\appr) \equiv 0$, thus completing the
proof of ii).
\par
It remains to consider the case where $d_\ast=\dtrig$. By the
arguments from the previous paragraph, it suffices to show that
$\dot{d}(t_\ast+\appr) \geq 0$ and $d(t_\ast+\appr) \geq \dtrig $.
Suppose that $d(t_\ast+\appr) \not\geq \dtrig$, i.e., there exists a
sequence $\{t_i\}$ such that $t_i> t_\ast, d(t_i) < \dtrig \;
\forall i$ and $t_i \to t_\ast$ as $i \to \infty$. Since $d(t_\ast)
= \dtrig$, a proper decrease of every $t_i$ gives $\dot{d}(t_i)<0$
in addition. By \eqref{c.a}, \eqref{ddot.plus}, the inequality
$\dot{d}(t)<0$ is maintained and enhanced as $t$ decreases from
$t_i$, remaining in the domain $\{t: d(t) < \dtrig\}$. Since
$\dot{d}(t_\ast)=0$, there is $\tau_i \in (t_\ast,t_i)$ such that
$d(\tau_i) = \dtrig$ and $\dot{d}(t) <0 \; \forall t \in
[\tau_i,t_i)$. Hence $d(\tau_i - \appr) > \dtrig$ and if $i$ is
large enough, there is $\theta_i > t_i$ such that $d(\theta_i) =
\dtrig$ and $d(t) < \dtrig \; \forall t \in (\tau_i, \theta_i)$.
Furthermore, there is $s_i \in (\tau_i, \theta_i)$ such that
$\dot{d}(t) <0 \; \forall t \in (\tau_i,s_i), \dot{d}(s_i) =0$,
$\dot{d}(t) \geq 0 \; \forall t \in [s_i,\theta_i]$. Then
$\beta(\theta_i)
>0$ by \eqref{beta.1}.
We note that $\beta(t_\ast)=0 \Rightarrow
\zeta_{\mathscr{P}}(t_\ast)=0$ for the vehicle path $\mathscr{P}$
and so $\zeta_{\mathscr{P}}(t) \to 0$ as $t \to t_\ast$. This and
\eqref{geom.eq} (applied to $\mathscr{P}$) imply that the sign of
$\dot{\beta}$ is determined by the sign of the path curvature:
\begin{equation}
\label{curv.turn} u = \pm \ov{u} \Rightarrow \pm \dot{\beta} <0
\qquad \forall t \approx t_\ast.
\end{equation}
\par
Suppose that $\exists \tau_\ast \in [\tau_i, s_i) : \beta(\tau_\ast)
\geq 0$. Since $u(t) = - \ov{u}\; \forall t \in (\tau_i,s_i)$, we
see that $\beta(s_i)>0, \dot{d}(s_i)=0, d_s:=d(s_i)$. By
Lemma~\ref{lem.slide}, sliding motion along the $d_s$-equidistant
curve is commenced at $t=s_i$ and maintained while $\beta>0$,
whereas $\beta>0$ until $\theta_i$ (if $i$ is large enough) due to
\eqref{beta.1}. However, this is impossible since $d_s< \dtrig$ and
$d(\theta_i) = \dtrig$. This contradiction proves that $\beta(t)
<0\; \forall t \in [\tau_i, s_i)$. The same argument and the
established validity of ii.2) for $d_\ast:=d_s< \dtrig$ show that
$\beta(s_i)<0$. Since $\beta(\theta_i) >0$, there exists $c_i \in
(s_i,\theta_i)$ such that $\beta(c_i) =0$ and $\beta(t)>0\; \forall
t \in (c_i,\theta_i]$. If $\dot{d}(c)=0$ for some $c \in
(c_i,\theta_i)$, Lemma~\ref{lem.slide} assures that sliding motion
along the $d(c)$-equidistant curve is started at $t=c$ and is not
terminated until $t=\theta_i$, in violation of $d(\theta)=\dtrig$.
For any $t \in (c_i,\theta_i)$, we thus have $\dot{d}(t)>0$. Hence
$u(t) = \ov{u}$ by \eqref{c.a}, $\dot{\beta}<0$ by
\eqref{curv.turn}, and so $\beta(c_i)=0 \Rightarrow \beta(\theta_i)
< 0$, in violation of the above inequality $\beta(\theta_i)>0$. This
contradiction proves that $d(t_\ast+\appr) \geq \dtrig $.
\par
Now suppose that $\dot{d}(t_\ast+\appr) \not\geq 0 $. Then there is
a sequence $\{t_i\}$ such that $t_i> t_\ast, \dot{d}(t_i) >0 \;
\forall i$ and $t_i \to t_\ast$ as $i \to \infty$; a proper increase
of every $t_i$ gives $d(t_i)> \dtrig$ in addition. By
\eqref{beta.1}, $d(t)> \dtrig \wedge \dot{d}(t)>0 \Rightarrow
\beta(t)>0$ for $t \approx t_\ast$ and so $u(t) = \ov{u}$ by
\eqref{c.a} and $\ddot{d}(t)<0$ by \eqref{ddot.plus}. So as $t$
decreases from $t_i$ to $t_\ast$, the derivative $\dot{d}(t)>0$
increases while $d> \dtrig$, in violation of the implication $d(t) =
\dtrig \Rightarrow \dot{d}(t) =0$ for $t \in [t_\ast,t_i]$. This
contradiction completes the proof.
\par
{\bf iii)} Were there a sequence $\{t_i\}_{i=1}^\infty$ such that
$\dot{d}(t_i)>0, \beta(t_i)>0 \; \forall i$ and $t_i \to t_\ast+0$
as $i \to \infty$, \eqref{c.a}, \eqref{ddot.plus}, and
\eqref{curv.turn} imply that as $t$ decreases from $t_i$ to $t_\ast$
for large enough $i$, the inequalities $\dot{d}(t)>0, \beta(t)>0$
would be preserved, in violation of $\dot{d}(t_\ast)=0,
\beta(t_\ast)=0$. It follows that $\dot{d}(t)>0 \Rightarrow \beta(t)
\leq 0$ for $t \approx t_\ast, t>t_\ast$.
\par Now assume existence of the sequence such that
$\dot{d}(t_i)>0, \beta(t_i) \leq 0 \; \forall i$ and $t_i \to
t_\ast+0$ as $i \to \infty$. For large $i$ such that $\beta(t_i)<0$,
\eqref{c.a}$\wedge$\eqref{ddot.plus} $\Rightarrow u(t) = - \ov{u}$,
and $\dot{d}(t)$ increases and so remains positive as $t$ grows from
$t_i$ until $\beta=0$. By \eqref{curv.turn},
$\ov{u}^{-1}|\beta(t_i)|$ time units later the vehicle becomes
headed to the target, which is trivially true if $\beta(t_i)=0$.
This and (i) of Lemma~\ref{lem.slide} imply that then the sliding
motion along the surface $\beta=0$ is commenced. It is maintained
while $\varkappa[s(t)] \geq 0$. Since $t_i \to t_\ast$ and
$\beta(t_i) \to \beta(t_\ast)=0$ as $i \to \infty$, this motion
occurs for $t> t_\ast$, i.e., iii) does hold.
\par
It remains to examine the case where $\dot{d}(t_\ast+\appr)\leq 0$
and so $d(t_\ast+\appr) \leq d_\ast$. Suppose first that either
$\dot{d}(t_\ast+\appr) \not\equiv 0$ or $\varkappa[s(t_\ast)+\appr]
\not\equiv 0$. Then $\beta(t_\ast+\appr) <0$ by \eqref{beta.2} and
$u = - \ov{u}$ at any side of the discontinuity surface $\dot{d}=0$
by \eqref{c.a}. Hence $u(t_\ast+\appr) \equiv - \ov{u}$, which
yields $\dot{d}(t_\ast+0)
>0$ by \eqref{ddot.plus}, in violation of $\dot{d}(t_\ast+0) = 0$. This
contradiction proves that $\dot{d}(t_\ast+\appr) \equiv 0$,
$\varkappa[s(t_\ast)+\appr]\equiv 0$. Then SMEC and SMT are
initially the same, and iii) does hold. \epf
\begin{remark}
\rm \label{rem.2} The times of switches between the modes of the discontinuous control law \eqref{c.a} do not accumulate.
\end{remark}
To prove this, we first note that the projection of any vehicle position
$\bldr$ within mode $\mathfrak{B}$ onto $\partial D$ is well defined
due to \eqref{geom.eq}. Let $s^-_i$ and $s^+_i$ be its values at the
start and end of the $i$th occurrence of the mode, respectively. By
Lemma~\ref{lem.sucs} and (vi) of Proposition~\ref{lem.cir}, $s$
monotonically sweeps an arc $\gamma_i$ of $\partial D$ with the ends
$s^-_i, s^+_i$ during the concluding part of $\mathfrak{B}$.
\begin{definition}
\label{def.single} The vehicle path or its part is said to be {\em
single} if the interiors of the involved arcs $\gamma_i$ are
pairwise disjoint and in the case of only one arc, do not cover
$\partial D$.
\end{definition}
Let $P$ and $Q$ be the numbers of the connected components of $S_\varkappa:= \{s : \varkappa(s) <0 \}$ and
$S_\zeta:=\{s : \zeta_{\partial D}(s) = 0 \}$, respectively. They are finite due to Corollary~\ref{cor.finite}.
\begin{lem}
\label{lem.ffnn} Any single path accommodates no more than
$(P+1)(Q+2)$ SMT's.
\end{lem}
\pf As was shown in the proof of (v) in of
Proposition~\ref{lem.cir}, the number of SMT's within a common mode
$\mathfrak{B}$ does not exceed $P+1$. SMT between the $i$th and
$(i+1)$th occurrences of $\mathfrak{B}$ starts at a position
$s_\dagger \in \gamma_i= [s^-_i, s^+_i]$ where $\zeta_{\partial
D}(s_\dagger) = - d <0$ and ends at the position $s^-_{i+1}$ where
$\zeta_{\partial D}(s^-_{i+1}) \geq 0$. Hence any arc $\gamma_i$,
except for the first and last ones, intersects adjacent connected
components $\text{Cc}_i^{=}$ and $\text{Cc}_i^{<}$ of $S_\zeta$ and
$\{s: \zeta_{\partial D}(s) <0\}$, respectively, such that the left
end-point of $\text{Cc}_i^{=}$ is the right end-point of
$\text{Cc}_i^{<}$. Hence $\text{Cc}_i^{=} \neq
\text{Cc}_{i^\prime}^{=} \; \forall i \neq i^\prime$, and so the
total number of the arcs $\gamma_i$ does not exceed $Q+2$, which
competes the proof. \epf
\par {\it Proof of Remark~\ref{rem.2}.}
Suppose to the contrary that the times $t_i$ when $\sigma$ is
updated accumulate, i.e., $t_i<t_{i+1} \to t_\ast < \infty$ as $i
\to \infty$. At $t=t_i$, a SMT is terminated, and so $d(t_i) =
\dtrig, \dot{d}(t_i) \leq 0, \beta(t_i) =0$. During the subsequent
AM, $d \leq \dtrig$. At such distances, \eqref{eq1} implies that
$|\ddot{d}| \leq M_d, |\ddot{s}| \leq M_s$, where $M_d, M_s >0$ do
not depend on the system state. Since IT ends with $\dot{d}=0$, this
AM lasts no less than $M_d^{-1}|\dot{d}(t_i)|$ time units. Hence
$\dot{d}(t_i) \to 0$ as $i \to \infty$. This and \eqref{eq1} imply
that $\dot{s}(t_i) - v \sgn \dot{s}(t_i) \to 0$ as $i \to \infty$.
So far as IT lasts no less than $M_s^{-1} |\dot{s}(t_i)|$ time units
if $\dot{s}$ is reversed during IT, the sign of $\dot{s}(t)$ is the
same for $t_i <t<t_\ast$ and large enough $i$. So the related part
of the path is single. By Lemma~\ref{lem.ffnn}, this part can
accommodate only a finite number of SMT's, in violation of the
initial hypothesis. This contradiction completes the proof.
\section{Proof of (ii) in Theorem~\ref{th.main}} \label{sec.prth}
This claim is identical to {\bf Remark~4$^{\boldsymbol{a}}$} from the basic paper.
\par
We first alter the control strategy by
replacement of the random machinery of choosing the turn direction
$\sigma$ at switches $\mathfrak{A} \mapsto \mathfrak{B}$ by a
deterministic rule. Then we show that the altered strategy achieves
the control objective by making no more than $N$ switches, where $N$
does not depend on the initial state of the robot. However, this
strategy cannot be implemented since it uses unavailable data. The
proof is completed by showing that with probability $1$, the initial
randomized control law sooner or later gives rise to $N$ successive
switches identical to those generated by the altered strategy.
\subsection{Deterministic Algorithm and its Properties}
\label{subsec.cl} \noindent The symbol $[\bldr_1,\bldr_2]$ stands
for the straight line segment directed from $\bldr_1$ to $\bldr_2$;
$\gamma_1 \star \gamma_2$ is the concatenation of directed curves
$\gamma_1, \gamma_2$ such that $\gamma_1$ ends at the origin of
$\gamma_2$.
\par
Let an occurrence $\mathfrak{A}^\dagger$ of mode $\mathfrak{A}$
holds between two modes $\mathfrak{B}$ and let it start at
$\bldr_\lozenge =\bldr(t_\lozenge)$ and end at $\bldr_\ast =
\bldr(t_\ast)$. Due to \eqref{eq.trig}, $\dist{\bldr_\ast} =
\dist{\bldr_\lozenge} = \dtrig$ are attained at unique boundary
points $s_\lozenge$ and $s_\ast$, respectively. They divide $C$ into
two arcs. Being concatenated with $\eta:= [s_\ast, \bldr_\ast] \star
[\bldr_\ast, \bldr_\lozenge] \star [\bldr_\lozenge, s_\lozenge]$,
each of them gives rise to a Jordan curve encircling a bounded
domain, one of which is the other united with $D$. The smaller
domain is denoted $\mathfrak{C}_{\mathfrak{A}^\dagger}$; it is
bounded by $\eta$ and one of the above arcs
$\gamma_{\mathfrak{A}^\dagger}$. Let $\sigma_{\mathfrak{A}^\dagger}
= \pm$ be the direction (on $\partial D$) of the walk from
$s_\lozenge$ to $s_\ast$ along $\gamma_{\mathfrak{A}^\dagger}$.
\par
We introduce the control law $\mathscr{A}$ that is the replica of
\eqref{c.a} except for the rule to update $\sigma$ when
$\mathfrak{A} \mapsto \mathfrak{B}$. Now for the first such switch,
$\sigma$ is set to an arbitrarily pre-specified value.
After any subsequent occurrence $\mathfrak{A}^\dagger$ of this mode,
\begin{equation}
\label{sigma.update} \sigma:= \left\{
\begin{array}{r l}
\sigma_{\mathfrak{A}^\dagger} &\text{if
$\mathfrak{C}_{\mathfrak{A}^\dagger}$ does not contain the target}
\\
-\sigma_{\mathfrak{A}^\dagger} & \text{if
$\mathfrak{C}_{\mathfrak{A}^\dagger}$ contains the target}
\end{array} \right. .
\end{equation}
\begin{prop}
\label{prop.det} Under the law $\mathscr{A}$, the target is reached
for a finite time, with making no more than $N$ switches
$\mathfrak{A} \mapsto \mathfrak{B}$, where $N$ does not depend on
the vehicle initial state.
\end{prop}
\par
The next two subsections are devoted to the proof of
Proposition~\ref{prop.det}. In doing so, the idea to retrace the
arguments justifying global convergence of the algorithms like the
Pledge one \cite{AbSe80} that deal with unconstrained motion of an
abstract point
is troubled by two problems. Firstly, this idea assumes that
analysis can be boiled down to study of a point moving according to
self-contained rules coherent in nature with the above algorithms.
i.e., those like 'move along the boundary', 'when hitting the
boundary, turn left', etc. However, this is hardly possible, at
least in full, since the vehicle behavior essentially depends on its
distance from the boundary. For example, depending on this distance
at the end of mode $\mathfrak{B}$, the vehicle afterwards may or may
not collide with a forward-horizon cusp of the obstacle. Secondly,
the Pledge algorithm and the likes are maze-escaping strategies;
they do not find the target inside a labyrinth when started outside
it. Novel arguments and techniques are required to justify the
success of the proposed algorithm in this situation.
\par
In what follows, we only partly reduce analysis of the vehicle
motion to that of a kinematically controlled abstract point. This
reduction concerns only special parts of the vehicle path and is not
extended on the entire trajectory. The obstacle to be avoided by the
point is introduced a posteriori with regard to the distance of the
real path from the real obstacle. To justify the convergence of the
abstract point to the target, we develop a novel technique based on
induction argument.
\par
We start with study of kinematically controlled point.
 \subsection{The Symbolic Path and its Properties}
 \label{subsect.spath}
 In this subsection, 'ray' means 'ray emitted from the target',
and we consider a domain $\mathscr{D}$ satisfying the following.
\begin{assumption}
\label{ass.finite1}
 The boundary $C:= \partial \mathscr{D}$ consists of finitely many
 (maybe, zero) straight line segments and the remainder on which the
 curvature vanishes no more than finitely many times. The domain
 $\mathscr{D}$ does not contain the target.
\end{assumption}
We also consider a point $\bldr$ moving in the plane according to
the following rules:
\begin{enumerate}[{\bf r.1)}]
\item The point moves outside the interior of $\mathscr{D}$;
\item Whenever $\bldr \not\in \mathscr{D}$, it moves
to $\bt$ in a straight line; 
\item Whenever $\bldr$ hits $\partial \mathscr{D}$,
it proceeds with monotonic motion along the boundary, counting the
angle $\beta$; \label{rule.smt}
\item This motion lasts until $\beta=0$ and new SMT is
possible, then SMT is commenced;
\item The point halts as soon as it arrives at the target.
\end{enumerate}
The possibility from r.4) means that $\mathscr{D}$ does not obstruct
the initial part of SMT. When passing the corner points of $\partial
\mathscr{D}$, the count of $\beta$ obeys \eqref{eq.angle} and the
conventional rules adopted for turning angles of the tangential
vector fields \cite{Krey91}, and is assumed to instantaneously,
continuously, and monotonically run between the one-sided limit
values. The possibility from r.4) may appear within this interval.
\par
To specify the turn direction in r.\ref{rule.smt}), we need some
constructions. Let the points $s_\pm \in C$ lie on a common ray and
$(s_-,s_+) \cap C = \emptyset$. One of them, say $s_-$, is closer to
the target than the other. They divide $C$ into two arcs. Being
concatenated with $(s_-,s_+)$, each arc gives rise to a Jordan curve
encircling a bounded domain. One of these domains is the other
united with $D$. The smaller domain $\mathfrak{C}(s_-,s_+)$ is
called the {\it cave} with the {\it corners} $s_-,s_+$. It is
bounded by $(s_-,s_+)$ and one of the above arcs
$\gamma_{\mathfrak{C}}$.
\par
To complete the rule r.\ref{rule.smt}), we note that any SMT except
for the first one starts and ends at some points $s_\lozenge, s_\ast
\in C$, which cut out a cave $\mathfrak{C}[s_\lozenge, s_\ast]$.
\begin{itemize}
\item[{\bf r.3a)}] After the first SMT, the turn is
in an arbitrarily pre-specified direction;
\item[{\bf r.3b)}] After SMT that is not the first the point turns
\begin{itemize}
\item outside $\mathfrak{C}[s_\lozenge, s_\ast]$ if the cave does
not contain the target;
\item inside the cave $\mathfrak{C}[s_\lozenge, s_\ast]$ if the cave
contains the target.
\end{itemize}
\end{itemize}
\begin{definition}
The path traced by the point obeying the rules {\rm r.1)---r.5),
r.3a), r.3b)} is called the {\em symbolic path (SP)}.
\end{definition}
\begin{prop}
\label{prop.sp} SP arrives at the target from any initial position.
The number of performed SMT's is upper limited by a constant $N$
independent of the initial position.
\end{prop}
The remainder of the subsection is devoted to the proof of this
claim. The notations $s,T,N, \fr, \lambda, \zeta$, $\varkappa$,
$\psi$, $\varphi$ are attributed to $C = \partial \mathscr{D}$. At
the corner points of $C$, these variables except for $s$ have
one-sided limits and are assumed to instantaneously, continuously,
and monotonically run between the one-sided limit values. An arc of
$C$ is said to be {\it regular} if $\zeta$ (non-strictly) does not
change its sign on this arc, depending on which the arc is said to
be {\it positive/negative} (or $\pm$arc). The regular arc is {\it
maximal} if it cannot be extended without violation of the
regularity. A connected part of $C$ and its points are said to be
{\it singular} if $\zeta$ strictly changes the sign when passing it
and, if this part contains more than one point, is identically zero
on it; see Fig.~\ref{fig.cross}(c). The singular arc is a segment of
a straight line since $\varkappa \equiv 0$ on it due to
\eqref{geom.eq}. The ends of any maximal regular arc are singular.
Due to Assumption~\ref{ass.finite1} and \eqref{geom.eq}, the
boundary $C$ has only finitely many singular parts. A boundary point
$s \in C$ is said to {\it lie above} $\mathscr{D}$ if there exists
$\delta>0$ such that $((1-\delta)s,s) \subset \mathscr{D}$ and
$(s,(1+\delta)s) \cap \mathscr{D} = \emptyset$. If conversely
$((1-\delta)s,s)\cap\mathscr{D} = \emptyset$ and $(s,(1+\delta)s)
\subset \mathscr{D}$, the point is said to {\it lie below}
$\mathscr{D}$.
\par
Formulas \eqref{geom.eq} and \eqref{angle.rot} imply the following.
\begin{observation}
\label{obs.arc} As $s$ moves in direction $\sigma=\pm$ over a
$\eta$-arc ($\eta=\pm$) of $C$, we have $\sigma \eta \dot{\varphi}
\geq 0$. Any point of $\pm$arc that is not singular lies above/below
$\mathscr{D}$.
\end{observation}
\begin{lem}
\label{lem.countbeta} As $s$ continuously moves along a regular arc,
$\beta$ evolves within an interval of the form $\Delta:= [\pi k, \pi
(k+1)]$, where $k$ is an integer. When $s$ reaches a singular point,
$\beta$ arrives at the end of $\Delta$ associated with the even or
odd integer, depending on whether $s$ moves towards or outwards the
target at this moment, respectively.
\end{lem}
\pf Since $\zeta$ does not change its sign, the vector $\fr$ does
not trespass the $\lambda$-axis, whereas $\beta$ is the polar angle
of this vector. This gives rise to the first claim of the lemma. The
second one is immediate from the first claim. \epf
\begin{lem}
\label{lem.beta} Whenever SP progresses along $C$ in direction
$\sigma = \pm$, we have $\sigma \beta \geq 0$.
\end{lem}
\pf This is evidently true just after any SMT. During the subsequent
motion along $C$, the inequality can be violated only at a position
$s$ where $\beta=0$ and either $s$ is a corner singular point or
$\varkappa(s+\sigma \appr) > 0$ since $\varkappa(s+\sigma \appr)
\leq 0 \Rightarrow \sigma \beta(s+\sigma \appr) \geq 0$ by the third
relation from \eqref{geom.eq}. However, at such position, motion
along $C$ is ended. \epf
\par
The cave $\mathfrak{C}(s_-,s_+)$ is said to be {\it
positive/negative} (or $\pm$cave) if the trip from $s_-$ to $s_+$
over $\gamma_{\mathfrak{C}}$ is in the respective direction of $C$.
By Observation~\ref{obs.arc}, $s$ moves from a $+$arc to a $-$arc in
this trip and so passes a singular part of $C$. The total number of
such parts inside $\gamma_{\mathfrak{C}}$ is called the {\it degree}
of the cave.\footnote{Possible singular parts at the ends of
$\gamma_{\mathfrak{C}}$ are not counted.}
\begin{lem}
For any cave of degree $M=1$, the arc $\gamma:=
\gamma_{\mathfrak{C}}$ consists of the positive $\gamma|_{s_- \to
s_-^\ast}$ and negative $\gamma|_{s_+^\ast \to s_+}$ sub-arcs and a
singular part $[s_-^\ast,s_+^\ast]$. For $s \in
[s_-^\ast,s_+^\ast]$, the tangential vector $T(s)$ (that is
co-linear with $[\bt,s]$ if $s$ is the corner point) is directed
outwards $\bt$ if the cave is positive and does not contain $\bt$ or
negative and contains $\bt$. Otherwise, this vector is directed
towards $\bt$. \label{lem.degone}
\end{lem}
\pf The first claim is evident. Let the cave be positive and $\bt
\not\in \mathfrak{C}(s_-,s_+)$. Suppose that $T(s)$ is directed
towards $\bt$. Then the same is true for $s:= s_+^\ast$. Hence
$\zeta(s_+^\ast+0) \leq 0$ and $\zeta(s_+^\ast+0) = 0 \Rightarrow
\lambda(s_+^\ast+0) >0 \Rightarrow \varkappa(s_+^\ast+\appr)
>0$ since otherwise, $\zeta (s_+^\ast+\appr) \geq 0$ by \eqref{geom.eq},
in violation of the definition of the singular part. In any case,
$((1-\delta)s_+^\ast, s_+^\ast ) \cap \mathscr{D} = \emptyset$ for
some $\delta >0$. Since $\bt \not\in \mathfrak{C}(s_-,s_+)$, the
segment $[0,s_+^\ast)$ intersects $\gamma_{\mathfrak{C}}$, cutting
out a smaller cave $\mathfrak{C}_{\text{sm}}$ inside
$\mathfrak{C}(s_-,s_+)$. The singular part inside
$\mathfrak{C}_{\text{sm}}$ is the second such part in the original
cave, in violation of $M=1$. This contradiction shows that $T(s)$ is
directed outwards $\bt$.
\par
Now suppose that $\bt \in \mathfrak{C}(s_-,s_+)$ and $T(s)$ is
directed outwards $\bt$. Let a point $s_\ast$ moves in the positive
direction along $\gamma|_{s_+^\ast \to s_+}$. The ray containing
$s_\ast$ monotonically rotates by Observation~\ref{obs.arc} and
contains a continuously moving point $s^{\text{mov}}_- \in
\gamma|_{s_-^\ast \to s_- }$. As $s_\ast$ runs from $s_+^\ast$ to
$s_+$, the segment $(s^{\text{mov}}_-, s_\ast)$ sweeps the entire
cave $\mathfrak{C}[s_-,s_+]$, and so this cave does not contain
$\bt$, in violation of the assumption. This contradiction proves
that $T(s)$ is directed towards $\bt$.
\par
The second claim for negative caves and the third claim are
established likewise. \epf
\begin{lem}
\label{lem.posbeta} If SP enters a cave without the target, it
leaves the cave through the other corner with $\beta \neq 0$. In
this maneuver, the direction of motion along $C$ is not changed, no
point of $C$ is passed twice, and the number of SMT's does not
exceed the cave degree.
\end{lem}
\pf Let SP enter the cave in the positive direction; the case of the
negative direction is considered likewise. The proof will be by
induction on the cave degree $M$.
\par
Let $M=1$. {\bf (i)} Suppose first that the cave is positive and so
$s$ enters it through $s_-$ moving over a $+$arc. By
Lemma~\ref{lem.degone}, the point $s$ moves outwards the target
whenever $s \in [s_-^\ast,s_+^\ast]$, and so $\beta \geq \pi$ by
Lemmas~\ref{lem.countbeta} and \ref{lem.beta}. As $s$ moves over the
subsequent $-$arc, $\zeta$ becomes negative and so the inequality is
kept true by Lemma~\ref{lem.countbeta}. Thus $s$ leaves the cave
through $s_+$ with $\beta \geq \pi
>0$, having made no SMT.
\par
{\bf (ii)} Let the cave be negative. Then $s$ enters it through
$s_+$ moving over the negative arc. By Lemma~\ref{lem.degone}, the
point $s$ moves towards the target whenever $s \in
[s_-^\ast,s_+^\ast]$. Since $\zeta(s_+ +0) \leq 0$,
Lemma~\ref{lem.beta} yields $\beta(s_++0) \geq \pi$. By
Lemma~\ref{lem.countbeta}, $\beta \geq \pi$ until $s_+^\ast$ and so
$\beta \geq 2 \pi$ at $s \in [s_-^\ast,s_+^\ast]$ by
Lemma~\ref{lem.degone}. When $s$ passes the entire
$[s_-^\ast,s_+^\ast]$, the sign of $\zeta$ reverses from $-$ to $+$
and so $\beta>2 \pi$ just after the passage of $s_-^\ast$. It
remains to note that $\beta \geq 2 \pi
>0$ while $s$ moves over the $+$arc from $s_-^\ast$ to $s_-$ by
Lemma~\ref{lem.countbeta}.
\par
Suppose that the claim of the lemma is true for any cave with degree
$\leq M$, and consider a cave of degree $M+1$. Let this cave be
positive. Then $s$ enters it through the lower corner $s_-$ along a
positive arc. We also consider the accompanying motion of the ray
containing $s$. This ray contains a continuously moving point
$s_+^\circledast \in C$ that starts at $s_+$. This motion is
considered until a singular part of $C$ appears on the ray segment
$[s,s_+^\circledast]$ for the first time. Three cases are possible
at this position.
\par
\begin{figure}
\subfigure[]{\scalebox{0.28}{\includegraphics{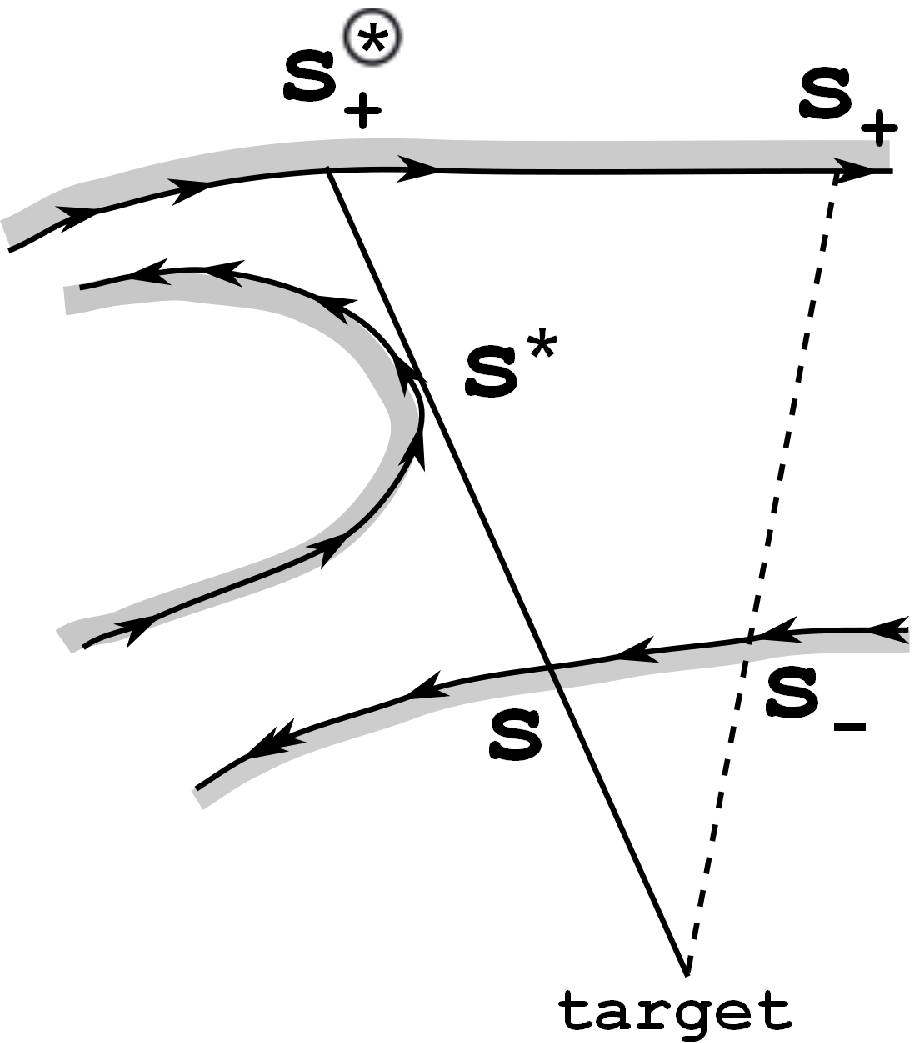}}}
\subfigure[]{\scalebox{0.28}{\includegraphics{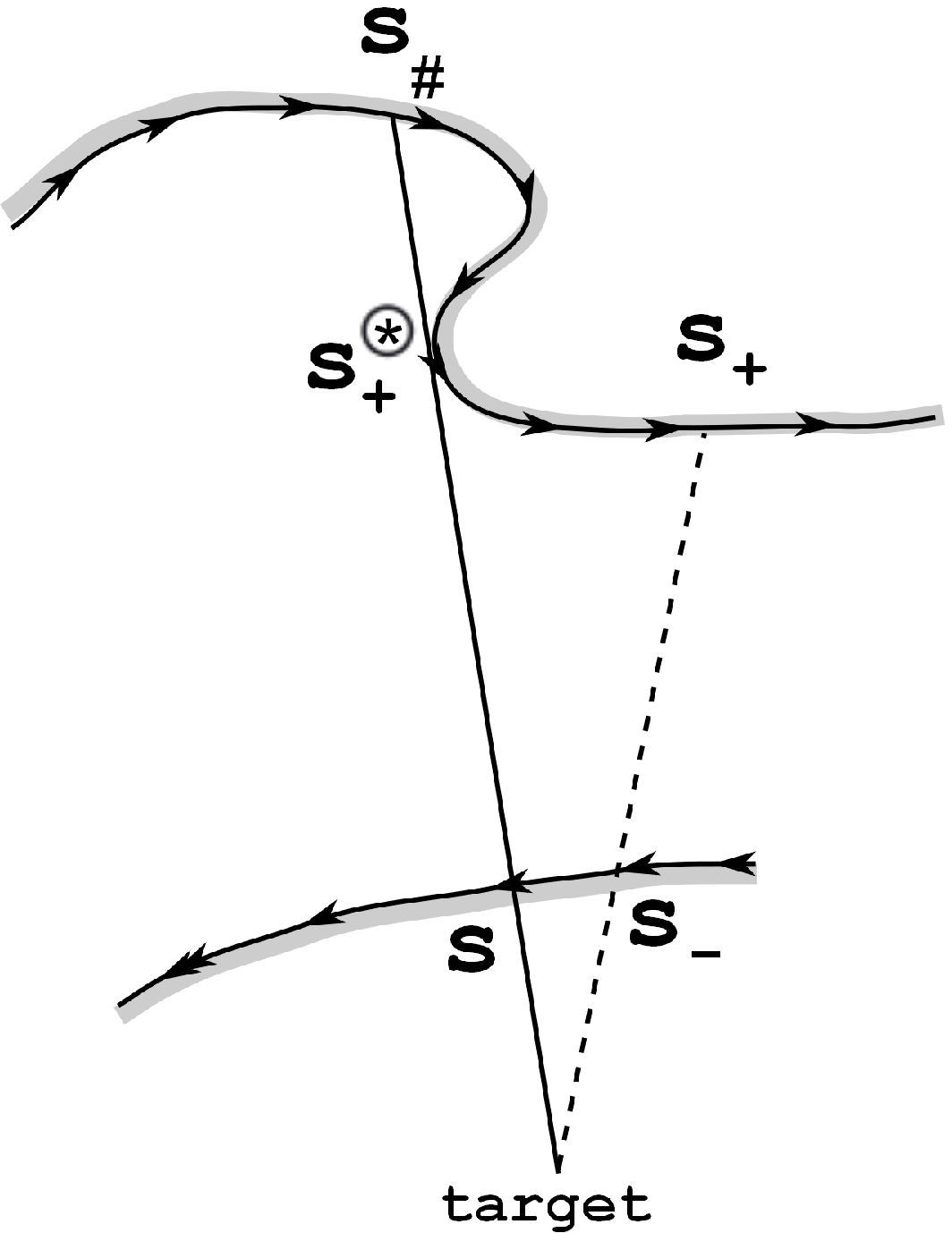}}}
\subfigure[]{\scalebox{0.28}{\includegraphics{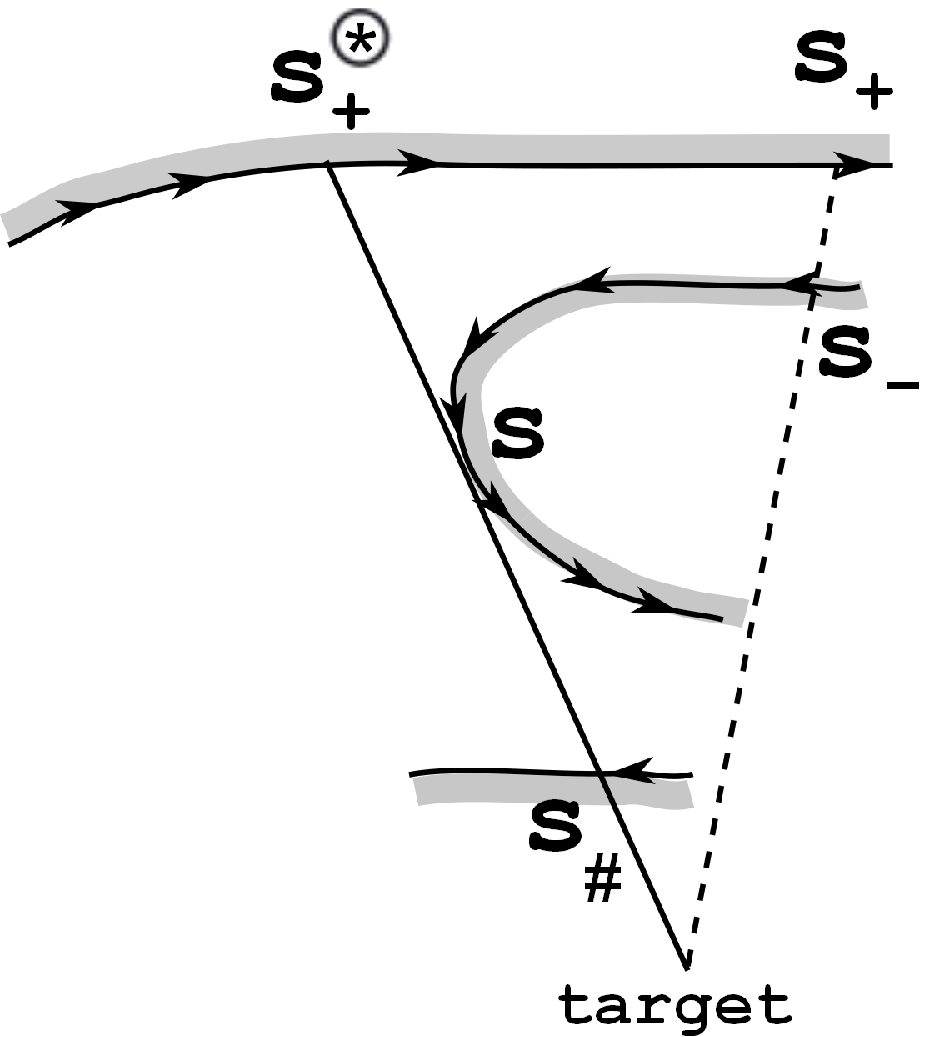}}}
\caption{The first singular point} \label{fig.frsing}
\end{figure}

{\bf (a)} The singular part $[s_-^\ast,s_+^\ast] \subset
(s,s_+^\circledast)$; see Fig.~\ref{fig.frsing}(a), where $s_-^\ast=
s_+^\ast=: s_\ast$.
By successively applying the induction hypothesis to
$\mathfrak{C}(s,s_+^\ast)$ and $\mathfrak{C}(s_+^\ast,
s^\circledast_+)$, we see that SP arrives at $s^\circledast_+$ in
the positive direction and with $\beta>0$. While $s$ moves from
$s^\ast_+$ to $s_+$ over the $-$arc, the vector $\fr(s)$ is below
the $\lambda$-axis and so $\beta \geq \pi >0$ by
Lemma~\ref{lem.countbeta}.
\par
{\bf (b)} The singular point $s_+^\circledast = s_-^\ast$; see
Fig.~\ref{fig.frsing}(b), where $s_-^\ast= s_+^\ast=: s_\ast$. By
successively applying the induction hypothesis to
$\mathfrak{C}(s,s_{\#})$ and $\mathfrak{C}(s_{\#},
s^\circledast_+)$, we see that SP arrives at $s^\circledast_+$ in
the positive direction and with $\beta>0$. So $\beta
(s^\circledast_+) \geq 2 \pi$ and SP proceeds along the $-$arc to
$s_+$ with $\beta \geq \pi >0$ by Lemma~\ref{lem.countbeta}, which
completes the proof.
\par
{\bf (c)} The singular point $s$; see Fig.~\ref{fig.frsing}(c). If
$\beta >0$ at this point, SP enters the cave
$\mathfrak{C}[s,s_+^\circledast]$ of degree $\leq M$ and by the
induction hypothesis, arrives at $s_+^\circledast$ moving in the
positive direction and with $\beta>0$. If conversely $\beta =0$, SP
undergoes SMT, which cannot be terminated at the target since it
does not belong to the cave at hand. So it is terminated at some
point $s_{\#} \in \gamma_{\mathfrak{C}}$. Since $\bt$ does not lie
in the sub-cave $\mathfrak{C}(s,s_{\#})$ of the original cave, the
vehicle turns right at $s_{\#}$ and thus proceeds along $C$ in the
positive direction. By applying the induction hypothesis to
$\mathfrak{C}(s_{\#}, s^\circledast_+)$, we see that SP arrives at
$s_+^\circledast$ moving in the positive direction and with
$\beta>0$ in any case. The proof is completed like in the cases (a)
and (b).
\par
The case where the cave is negative is considered likewise.
\begin{lem}
\label{lem.notchange} Suppose that after SMT starting and ending at
the points $s_\lozenge$ and $s_\ast$, respectively, the direction of
the motion along $C$ is reversed. Then the cave
$\mathfrak{C}[s_\lozenge,s_\ast]$ does not contain $\bt$ but
contains the entire path traced before SMT at hand.
\end{lem}
\pf Let the motion direction at $s=s_\lozenge$ be $+$; the case of
$-$ is considered likewise. Since on arrival at $s_\ast$, the left
turn is made,  $\mathfrak{C}[s_\lozenge,s_\ast]$ does not contain
$\bt$ by r.3b). Suppose that the path traced before SMT at hand is
not contained by this cave, i.e., the point enters this cave before.
Since this cannot be done during another SMT, the point enters the
cave through either $s_\lozenge$ or $s_\ast$. In the first case,
$s_\lozenge$ is passed twice in the opposite directions, in
violation of Lemma~\ref{lem.posbeta}. In the second case,
$s_\lozenge$ is passed with $\beta>0$ by the same lemma and so SMT
cannot be commenced. The contradiction obtained proves that the
initial part of SP is inside the cave. \epf
\begin{lem}
\label{lem.posbeta1} If SP progresses along $C$ in a cave not
containing the target, it leaves this cave through one of its
corners. During this maneuver, SP passes no point of $C$ twice and
makes no more SMT's than the degree of the cave.
\end{lem}
\pf For the definiteness, let the cave be positive; the case of the
negative cave is considered likewise. The proof will be by induction
on the degree $M$ of the cave.
\par
Let $M=1$. We employ the notations from Lemma~\ref{lem.degone}.
\par
{\bf ($\boldsymbol{\alpha}$)} The motion is started on
$\gamma|_{s_+^\ast \to s_-}$ in the direction $-$. The claim is
evident.
\par
{\bf ($\boldsymbol{\beta}$)} The motion is started on $\gamma|_{s_+
\to s_+^\ast}$ in the direction $-$. Then the point necessarily
arrives at $s_+^\ast$, moving in the negative direction. Thus the
situation is reduced to ($\alpha$).
\par
{\bf ($\boldsymbol{\gamma}$)} The motion is started on
$\gamma|_{s_-^\ast \to s_+}$ in the positive direction. The claim of
the lemma is justified by the concluding arguments from (i) in the
proof of Lemma~\ref{lem.posbeta}.
\par
{\bf ($\boldsymbol{\delta}$)} The motion is started on $\gamma|_{s_-
\to s_-^\ast}$ in the direction $+$. Then the point necessarily
arrives at $s_-^\ast$, moving in the positive direction. Thus the
situation is reduced to ($\gamma$).
\par
Now suppose that the claim of the lemma is true for any cave with
degree $\leq M$, and consider a cave of degree $M+1$. Let this cave
be positive for the definiteness; the case of the negative cave is
considered likewise. We also consider an auxiliary motion of the
point over $C$ from $s_-$ into the cave and the accompanying motion
of the ray containing $s$ until one of the situations from
Fig.~\ref{fig.frsing} occurs.
\par
{\bf Case (a) from Fig.~\ref{fig.frsing}.} {\bf (a.1)} If the motion
is started on $\gamma|_{s^\circledast_+ \to s_+}$ in direction $+$
or on $\gamma|_{s \to s_-}$ in direction $-$, the claim of the lemma
is justified by the concluding arguments from (i) in the proof of
Lemma~\ref{lem.posbeta}.
\par
{\bf (a.2)} If the motion is started on $\gamma|_{s_-^\ast \to
s^\circledast_+}$, the induction hypothesis applied to the cave
$\mathfrak{C}[s_-^\ast , s^\circledast_+]$ of degree $\leq M$
ensures that the point arrives at either $s^\circledast_+$ or
$s_-^\ast$. In the first case, it arrives in direction $+$, thus
reducing the situation to (a.1). In the second case, it arrives in
direction $-$. If $\beta \neq 0$ at this position, the point enters
the cave $\mathfrak{C}[s_-^\ast,s]$ in direction $-$ and afterwards
leaves it through $s$ in the same direction by
Lemma~\ref{lem.posbeta}. If $\beta=0$, SMT is commenced, which ends
at the position $s$ with the left turn since
$\mathfrak{C}[s_-^\ast,s]$ does not contain $\bt$. Hence in any
case, the motion proceeds in direction $-$ from the position $s$,
which reduces the situation to (a.1).
\par
{\bf (a.3)} The case where the motion is started on $\gamma|_{s \to
s_-^\ast}$, is considered likewise.
\par
{\bf (a.4)} The cases where the motion starts on
$\gamma|_{s^\circledast_+ \to s_+}$ in direction $-$ or on
$\gamma|_{s \to s_-}$ in direction $+$, are trivially reduced to
(a.2) and (a.3), respectively.
\par
{\bf Case (b) from Fig.~\ref{fig.frsing}.} {\bf (b.1)} The cases
where the motion starts on $\gamma|_{s^\circledast_+ \to s_+}$ in
direction $+$ or on $\gamma|_{s \to s_-}$ in direction $-$, is
considered like (a.1).
\par
{\bf (b.2)} If the start is on $\gamma|_{s \to s_{\#}}$, the
induction hypothesis applied to $\mathfrak{C}[s , s_{\#}]$ ensures
that the point arrives at either $s$ or $s_{\#}$. In the first case,
it arrives in direction $-$, thus reducing the situation to (b.1).
In the second case, it arrives in direction $+$ and then enters the
cave $\mathfrak{C}[s_{\#}, s_+^\circledast]$. By
Lemma~\ref{lem.posbeta}, the point leaves this cave through
$s_+^\circledast$ in direction $+$ and with $\beta >0$, thus
reducing the situation to (b.1).
\par
{\bf (b.3)} If the motion commences on $\gamma|_{s_{\#} \to
s^\circledast_+}$, the induction hypothesis applied to the cave
$\mathfrak{C}[s_{\#}, s^\circledast_+]$ of degree $\leq M$ ensures
that the point arrives at either $s_{\#}$ or $s^\circledast_+$. In
the first case, the arrival is in direction $-$, after which the
situation is reduced to (b.2). In the second case, the arrival is in
direction $+$. If $\beta \neq 0$ at this moment, the motion proceeds
along $\gamma|_{s^\circledast_+ \to s_+}$ in direction $+$, and the
situation is reduced to (b.1). If $\beta =0$, SMT is commenced,
which ends at the position $s$ with the left turn since the cave
$\mathfrak{C}[s^\circledast_+,s]$ does not contain the target. Hence
the motion proceeds along $\gamma|_{s \to s_-}$ in direction $-$,
and the situation is still reduced to (b.1).
\par
{\bf (b.4)} The cases where the motion starts on
$\gamma|_{s^\circledast_+ \to s_+}$ in direction $-$ or on
$\gamma|_{s \to s_-}$ in direction $+$, are trivially reduced to
(b.3) and (b.2), respectively.
\par
{\bf Case (c) from Fig.~\ref{fig.frsing}.} {\bf (c.1)} The cases
where the motion starts on $\gamma|_{s^\circledast_+ \to s_+}$ in
direction $+$ or on $\gamma|_{s \to s_-}$ in direction $-$, is
considered like (a.1).
\par
{\bf (c.2)} If the start is on $\gamma|_{s_{\#} \to
s_+^\circledast}$, the induction hypothesis applied to
$\mathfrak{C}[s_{\#}, s_+^\circledast]$ yields that the point
arrives at either $s_+^\circledast$ or $s_{\#}$. In the first case,
the arrival direction is $+$ and the situation is reduced to (b.1).
In the second case, the point arrives in direction $-$ and then
enters $\mathfrak{C}[s_{\#}, s]$. By Lemma~\ref{lem.posbeta}, the
point leaves this cave through $s$ in direction $-$ and with $\beta
>0$. Thus we arrive at (b.1) once more.
\par
{\bf (c.3)} If the motion commences on $\gamma|_{s_{\#} \to s}$, the
induction hypothesis applied to the cave $\mathfrak{C}[s_{\#}, s]$
of degree $\leq M$ ensures that the point arrives at either $s_{\#}$
or $s$. In the first case, the arrival is in direction $+$, after
which the situation is reduced to (b.2). In the second case, the
arrival is in direction $-$, after which the situation reduces to
(b.1).
\par
{\bf (c.4)} The cases where the motion starts on
$\gamma|_{s^\circledast_+ \to s_+}$ in direction $-$ or on
$\gamma|_{s \to s_-}$ in direction $+$, are trivially reduced to
(c.2) and (c.3), respectively. \epf
\begin{lem}
\label{lem.end} Any part of SP where it progresses over the boundary
$\partial D$ ends with SMT.
\end{lem}
\pf is by retracing the proof of (v) in Proposition~\ref{lem.cir}.
\par
Let $K$ be the number of singular parts of the boundary $\partial
\mathscr{D}$.
\begin{lem}
\label{lem.ncont0} If every cave examined in {\rm r.3b)} does not
contain the target, SP consists of the initial $\mathscr{P}^-$ and
terminal $\mathscr{P}^+$ sub-paths (some of which may contain only
one point) such that each accommodates no more than $K$ SMT's, no
point of $C$ is passed twice within $\mathscr{P}^-$, whereas the
direction of motion along $C$ is not altered within $\mathscr{P}^+$.
\end{lem}
\pf Suppose first that the initial position lies in some cave. Among
such caves, there is one enveloping the others. By
Lemma~\ref{lem.posbeta1}, SP leaves this cave and the related
sub-path satisfies the properties stated in Lemma~\ref{lem.ncont0}.
If the initial position lies outside any cave, this sub-path is
taken to consist of only this position. By
Lemma~\ref{lem.notchange}, the direction of the motion along $C$ is
not changed on the remaining sub-path $\mathscr{P}_+$ and
$\mathscr{P}_+$ does not go inside the above maximal cave.
\par
Suppose that within $\mathscr{P}_+$, SP accommodates more than $K$
SMT's. Any of them starts at some singular part with $\beta=0$.
Hence SP passes some singular point with $\beta=0$ at least twice
and thus becomes cyclic. Now we consider the related minimal cyclic
part CP of SP that starts and ends with commencing a SMT at a common
point. Due to the constant direction, the closed curve CP is simple.
It follows that $\tang{\text{CP}} = \pm 2 \pi$, whereas
$\sphan{\bt}{\text{CP}} =0$ since $W=0$ for all bypassed caves and
$\bt \not\in D$. Hence $\sphan{0}{\fr}=\mp 2 \pi$ by
\eqref{eq.angle}, whereas CP starts and ends with $\beta=0$ and so
$\sphan{0}{\fr}= 0$. This contradiction completes the proof. \epf
\par
Lemmas~\ref{lem.end} and \ref{lem.ncont0} give rise to the
following.
\begin{corollary}
\label{lem.ncont} If every cave examined in {\rm r.3b)} does not
contain $\bt$, SP arrives at $\bt$ by making no more than $2K$
SMT's.
\end{corollary}
\begin{lem}
\label{lem.posbeta2} If SP enters a cave containing $\bt$ over a
positive arc with $|\beta| \leq \pi$, it arrives at $\bt$ not
leaving the cave. During this maneuver, no point of $C$ is passed
twice and the number of SMT's does not exceed the degree of the
cave.
\end{lem}
\pf Let the cave be entered in direction $+$; the case of $-$ is
considered likewise. The proof will be by induction on the degree
$M$ of the cave $\mathfrak{C}[s_-,s_+]$. Since $s$ enters the cave
over a positive arc, the entrance is through $s_-$.
\par
Let $M=1$. By Lemma~\ref{lem.degone}, $s$ moves towards $\bt$ when
reaching the singular part of the cave $[s_-^\ast, s_+^\ast]$. At
this position, $\beta=0$ by Lemma~\ref{lem.countbeta} and
$\mathscr{D}$ does not obstruct the initial part of SMT, as was show
in the proof of Lemma~\ref{lem.degone}. So SMT is commenced. If it
is not terminated at $\bt$, the segment $[0,s_-^\ast)$ intersects
$\gamma_{\mathfrak{C}}$, cutting out a smaller cave within the
original one. The singular part inside this new cave is the second
such part within the original cave, in violation of $M=1$. Hence
$\bt$ is reached and only one switch $\mathfrak{B} \mapsto
\mathfrak{A}$ is made.
\par
Now suppose that the conclusion of the lemma is true for any cave
with degree $\leq M$, and consider a cave of degree $M+1$. Like in
the proof of Lemma~\ref{lem.posbeta}, we consider the motion of the
ray containing $s$ until a singular point appears on the segment
$[s,s_+^\ast]$ for the first time, and examine separately three
possible cases depicted in Fig.~\ref{fig.frsing}.
\par
\par
{\bf (a)} The singular point $s_\ast \in (s,s_+^\ast)$; see
Fig.~\ref{fig.frsing}(a). The target is contained by the cave
$\mathfrak{C}[s,s_\ast]$ of degree $\leq M$, which is entered in the
positive direction and by Lemma~\ref{lem.countbeta}, with $0 \leq
\beta \leq \pi$. The induction hypothesis competes the proof.
\par
{\bf (b)} The singular point $s_\ast = s_+^\ast$; see
Fig.~\ref{fig.frsing}(b). The target is evidently contained by the
cave $\mathfrak{C}[s,s_{\#}]$ of degree $\leq M$. The proof is
completed like in the previous case.
\par
{\bf (c)} The singular point $s_\ast = s$; see
Fig.~\ref{fig.frsing}(c). If at $s_\ast$, the point moves outwards
$\bt$, the arguments from the second paragraph in the proof of
Lemma~\ref{lem.degone} show that the cave does not contain $\bt$, in
violation of the assumption of the lemma. Hence at $s_\ast$, the
point moves towards $\bt$ and so $\beta =0$ by
Lemma~\ref{lem.countbeta} and $\mathscr{D}$ does not obstruct the
initial part of SMT, as was show in the proof of
Lemma~\ref{lem.degone}. Thus SMT is commenced at $s_\ast$. If it is
terminated at $\bt$, the proof is completed. Otherwise, it arrives
at $s_{\#} \in \gamma_{\mathfrak{C}}$, as is shown in
Fig.~\ref{fig.frsing}(c). Evidently, the cave
$\mathfrak{C}[s_{\#},s]$ does not contain the target. So on reaching
$s_{\#}$, the point turns right and continues moving in the positive
direction over a new positive arc and with $\beta \in [0,\pi]$. So
the proof is completed by applying the induction hypothesis to the
cave $\mathfrak{C}[s_{\#}, s_+^\ast]$ of degree $\leq M$.
\par
{\it Proof of Proposition~\ref{prop.sp}:} is straightforward from
Corollary~\ref{lem.ncont} and Lemma~\ref{lem.posbeta2}.
\subsection{Proof of Proposition~\ref{prop.det}.}
Let $\mathscr{P}$ stand for the directed path traced by the vehicle
under the control law $\mathscr{A}$ from Subsect.~\ref{subsec.cl}.
We first show that after a slight modification, this path can be
viewed as SP for some domain $\mathscr{D}$ provided that
$\mathscr{P}$ is single (see Definition~\ref{def.single}). This
permits us to employ the results of Subsect.~\ref{subsect.spath}.
\par
We use the notations $s^-_i, s^+_i, \gamma_i$ from
introduced before Definition~\ref{def.single}, note that for $s \in \gamma_i$, the
distance $d$ from the vehicle to the obstacle is a function
$d=d_i(s)$ of $s$, and put:
\begin{multline}
\label{eq.d}
 \mathscr{D} := \Big\{ \bldr: d:=
\dist{\bldr} < d_\star(D)\; \text{an either} \;  s:=  \\
s(\bldr) \in \gamma_i \wedge d \leq d_i(s) \; \text{or}\; s \not\in
\cup_i \gamma_i \wedge d \leq \dtrig \Big\}.
\end{multline}
If $\sigma \dot{s}<0$ at the start of the $i$th mode $\mathfrak{B}$,
the abscissa $s^-_i$ is passed twice during IT by
Lemma~\ref{lem.sucs}. For every such $i$, the real path between
these two passages is replaced by the motion along the straight line
segment, which gives rise to the modified path $\mathscr{P}_\ast$.
\begin{observation}
Let the original path be single. The modified path
$\mathscr{P}_\ast$ is SP for $\mathscr{D}_\ast$.
\end{observation}
Indeed, this path can be viewed as a trace of a point obeying the
rules r.1)---r.5). To ensure r.3a), the direction should be
pre-specified to match that of $\mathscr{P}_\ast$. The property
r3.b) is satisfied due to \eqref{sigma.update} and the second
inequality from \eqref{far.enough}.
\begin{lem}
\label{lem.ns} For a single path, the set \eqref{eq.d} satisfies
Assumption~{\rm \ref{ass.finite1}} and its boundary has no more than
$N_s$ singular parts, where $N_s$ is completely determined by $D$
and $\bt$.
\end{lem}
\pf The last claim in Assumption~\ref{ass.finite1} holds by
\eqref{far.enough}, \eqref{eq.d}. The boundary $\partial
\mathscr{D}$ consists of parts traced during 1) SMT's, 2) SMEC's, 3)
arcs of circles traced during IT's, and 4) segments of normals to
$\partial D$ resulted from the path modification.
\par
Any part 1) clearly satisfies Assumption~\ref{ass.finite1} and is
either singular or does not contain singular points; their number
does not exceed $(P+1)(Q+1)$ by Lemma~\ref{lem.ffnn}.
\par
Since parts 2) are separated by SMT's, their number does not exceed
$(P+1)(Q+1)+1$. Any part 2) lies on a $d$-equidistant curve $C(d)$
with $d\leq \dtrig$. Due to \eqref{eqdc}, $\zeta_{C(d)}(s) =
\zeta_{\partial D}(s)+d$, Assumption~\ref{ass.finite1} holds since the boundary $\partial D$ is piece-wise analytical, and the singular parts of $C(d)$ are
the connected components of the set from Corollary~\ref{cor.finite}.
So type 2) arcs of $C$ accommodate no more than $F[(P+1)(Q+1)+1]$
singular parts.
\par
It remains to note that parts 3) and 4) do not contain singular
points since $\beta$ monotonically evolves from $0$ during IT's.
\begin{lem}
\label{lem.find} If the vehicle finds the target in
$\mathfrak{C}_{\mathfrak{A}^\dagger}$ after some occurrence
$\mathfrak{A}^\dagger$ of mode $\mathfrak{A}$, it arrives at the
target by making after this no more than $N_s$ switches
$\mathfrak{A} \mapsto \mathfrak{B}$.
\end{lem}
\pf Let us consider a part $\mathscr{P}$ of the path that starts in
mode $\mathfrak{B}$ preceding $\mathfrak{A}^\dagger$. Suppose first
that this part is not single and truncate it from the right, leaving
its maximal single sub-part $\mathscr{P}^\dagger$. The terminal
position of $\mathscr{P}^\dagger$ lies on a previously passed piece
of $\mathscr{P}^\dagger$. Let $\mathscr{D}^\dagger$ and
$\mathscr{P}^\dagger_\ast$ be the related domain \eqref{eq.d} and
modified path. Associated with $\mathfrak{C}_{\mathfrak{A}^\dagger}$
is a cave of $\mathscr{D}^\dagger$ into which
$\mathscr{P}^\dagger_\ast$ turns with $|\beta| \leq \pi$. By
Lemma~\ref{lem.posbeta2}, $\mathscr{P}^\dagger_\ast$ cannot arrive
at a previously passed point, in violation of the above property.
This contradiction proves that the entire path $\mathscr{P}$ is
single. Then Lemmas~\ref{lem.posbeta2} and \ref{lem.ns} guarantee
that $\mathscr{P}_\ast$ arrives at $\bt$ by making no more than
$N_s$ SMT's. It remains to note that $\mathscr{P}$ and
$\mathscr{P}_\ast$ arrive at $\bt$ only simultaneously, and each
occurrence of $\mathfrak{A}$ gives rise to a SMT in
$\mathscr{P}_\ast$.
\begin{lem}
\label{lem.direct} After no more than $N_s+1$ switches $\mathfrak{A}
\mapsto \mathfrak{B}$, the direction in which $s$ moves along
$\partial D$ within modes $\mathfrak{B}$ is not altered.
\end{lem}
\pf Consider an occurrence $\mathfrak{A}^\dagger$ of mode
$\mathfrak{A}$ after which the direction is altered and the path
$\mathscr{P}$ from the start of the entire motion until the end of
$\mathfrak{A}^\dagger$. Suppose that $\mathscr{P}$ is not single and
truncate it from the left, leaving the maximal single part
$\mathscr{P}^\dagger$. The starting point of $\mathscr{P}^\dagger$
is passed once more within $\mathscr{P}^\dagger$, both times in mode
$\mathfrak{B}$. So this double point is inherited by
$\mathscr{P}^\dagger_\ast$, where $\mathscr{D}^\dagger$ and
$\mathscr{P}^\dagger_\ast$ are the related domain \eqref{eq.d} and
modified path. Associated with $\mathfrak{C}_{\mathfrak{A}^\dagger}$
is a cave $\mathfrak{C}_{\mathscr{D}^\dagger}$ of
$\mathscr{D}^\dagger$; these two sets contain the target only
simultaneously due to \eqref{far.enough}. Hence $\mathscr{P}$ and
$\mathscr{P}^\dagger_\ast$ acquire a common turn direction at their
ends. So SP $\mathscr{P}^\dagger_\ast$ has converse directions of
motion along the boundary at the start and end of the last involved
SMT and by Lemmas~\ref{lem.notchange} and \ref{lem.posbeta1}, has no
double points. This contradiction proves that the entire
$\mathscr{P}^\dagger$ is single. Due to Lemma~\ref{lem.notchange},
the modified path $\mathscr{P}_\ast^\dagger$ lies in
$\mathfrak{C}_{\mathscr{D}^\dagger}$ and so involves no more than
$N_s$ SMT's thanks to Lemmas~\ref{lem.posbeta1} and \ref{lem.ns}. It
remains to note that each occurrence of $\mathfrak{A}$ gives rise to
a SMT in $\mathscr{P}_\ast$. \epf
\par
To prove Proposition~\ref{prop.det}, it in fact remains to show that
the vehicle cannot pass more than $N_s$ modes $\mathfrak{A}$ in a
row, constantly not finding the target in
$\mathfrak{C}_{\mathfrak{A}}$ and not changing the direction of the
motion along $\partial D$. The next lemma with corollaries serves
this proof. The symbol $\angle (a, b) \in (-\pi, \pi]$ stands for
the angle from the vector $a$ to $b$. Let the points $\bldr_i,
i=1,2$ on $\mathscr{P}$ be at the distance $\dist{\bldr_i}\leq
\dtrig$ and such that when traveling between them, the path does not
intersect itself and except for $\bldr_i$, has no points in common
with the normals $[\bldr_i, s_i]$, where $s_i:=s[\bldr_i]$. The
points $s_i$ split $\partial D$ into two curves. Being concatenated
with the above normals and $\mathscr{P}|_{\bldr_1 \to \bldr_2}$,
they give rise to Jordan loops, with one of them enveloping the
other. Let $\gamma_{\text{inner}}$ be the curve giving rise to the
inner loop $\text{\bf LOOP}$, and $\sigma = \pm$ be the direction
from $s_1$ to $s_2$ along $\gamma_{\text{inner}}$.
\begin{lem}
\label{lem.turn2} If $\text{\bf LOOP}$ does not encircle the target,
the following relation holds
\begin{multline}
\label{basic.eq} \sphan{0}{\fr_{\mathscr{P}}|_{\bldr_1 \to \bldr_2}}
= \sphan{0}{\fr_{\partial D}|_{s_1 \xrightarrow{\sigma} s_2}} +
\sphan{\bt}{[\bldr_1, s_1]} - \sphan{\bt}{[\bldr_2, s_2]}
\\
+ \angle \left[ \sigma T_{\partial D}(s_1),
T_{\mathscr{P}}(\bldr_1)\right] - \angle \left[ \sigma T_{\partial
D}(s_2), T_{\mathscr{P}}(\bldr_2)\right].
\end{multline}
\end{lem}
\pf Let $\sigma=+$; $\sigma=-$ is considered likewise. By applying
the Hopf's theorem to {\bf LOOP}, we see that
\begin{multline*}
\sphan{\bt}{[s_1, \bldr_1]}+ \sphan{\bt}{\mathscr{P}|_{\bldr_1 \to
\bldr_2}} + \sphan{\bt}{[\bldr_2, s_2]} - \sphan{\bt}{\partial
D_{s_1 \to s_2}} =0,
\\
\tang{\mathscr{P}|_{\bldr_1 \to \bldr_2}} = \tang{\partial D_{s_1
\to s_2}} \\
- \angle \left[ T_{\partial D}(s_1), T_{\mathscr{P}}(\bldr_1)\right]
+ \angle \left[ T_{\partial D}(s_2),
T_{\mathscr{P}}(\bldr_2)\right].
\end{multline*}
The proof is completed by the second formula in \eqref{eq.angle}.
\epf
\par
\begin{figure}
\scalebox{0.25}{\includegraphics{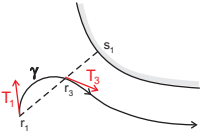}}
\caption{Auxiliary loop} \label{fig.loop}
\end{figure}

The next claim employs the notations introduced at the beginning of
Subsect.~\ref{subsec.cl}.
\begin{corollary}
\label{lem.turn1} Suppose that $\bt \not\in
\mathfrak{C}_{\mathfrak{A}^\dagger}$ and the value of $\sigma$
maintained during the occurrence $\mathfrak{A}^\dagger$ of mode
$\mathfrak{A}$ is not altered when $\mathfrak{A}^\dagger \mapsto
\mathfrak{B}$. Then \eqref{basic.eq} holds with $\bldr_1:=
\bldr_\lozenge, \bldr_2:= \bldr_\ast$.
\end{corollary}
This is true since in this claim and Lemma~\ref{lem.turn2}, $\sigma$
is the same.
\begin{corollary}
\label{lem.parr1} Let $\bldr_1$ and $\bldr_2$ be successively passed
within a common mode $\mathfrak{B}$, where $\sigma(t) \equiv
\sigma=\pm$. If $\bldr_2$ is passed after IT, \eqref{basic.eq}
holds, where $\sphan{0}{\fr_{\partial D}|_{s_1 \xrightarrow{\sigma}
s_2}}$ accounts for the entire motion of the projection $s=s[\bldr],
\bldr \in \mathscr{P}_{\bldr_1 \to \bldr_2}$, including possible
full runs over $\partial D$.
\end{corollary}
If 1) $s$ does not run over the entire $\partial D$ and 2) either
$\bldr_1$ is passed after IT or $\sgn \dot{s} = \sigma$ at the start
of the mode, the claim is evident. If 1) holds but 2) does not, the
path may intersect $[s_1,\bldr_1]$ and so direct application of
Lemma~\ref{lem.turn2} is impossible. Then we apply this lemma to
$\bldr_1:= \bldr_3$, where $\bldr_3$ is the point where the vehicle
intersects the normal for the second time during IT; see
Fig.~\ref{fig.loop}. The proof is completed by noting that
$\sphan{\bt}{[\bldr_1, \bldr_3]} = \sphan{\bt}{\gamma},
\tang{\gamma} = \angle[T_1,T_3]$ and so
$\sphan{0}{\fr_{\mathscr{P}}|_{\bldr_1 \to \bldr_3}} =
\sphan{0}{\fr_{\gamma}} = \sphan{\bt}{[\bldr_1, \bldr_3]}-
\angle[T_1,T_3]$, as well as that $\angle \left[ \sigma T_{\partial
D}(s_1), T_3\right] = \angle \left[ \sigma T_{\partial D}(s_1),
T_1\right]+ \angle[T_1,T_3]$. The claim is generalized on the case
where 1) is not true by proper partition of the path, followed by
summation of the formulas related to the resultant pieces.
\begin{corollary}
\label{cor.lasrt}
Let points $\bldr_1$ and $\bldr_2$ be successively passed in modes
$\mathfrak{B}$ (maybe, different). Suppose that $\bldr_2$ is not
attributed to IT and when traveling from $\bldr_1$ to $\bldr_2$, the
vehicle constantly does not find the target in
$\mathfrak{C}_{\mathfrak{A}}$ and does not change $\sigma$. Then
\eqref{basic.eq} holds, where $\sphan{0}{\fr_{\partial D}|_{s_1
\xrightarrow{\sigma} s_2}}$ accounts for the entire motion of the
projection $s=s[\bldr], \bldr \in \mathscr{P}_{\bldr_1 \to
\bldr_2}$, including possible full runs over $\partial D$.
\end{corollary}
It is assumed that as the vehicle moves in mode $\mathfrak{A}$, the
projection $s$ continuously and monotonically goes over $\partial D$
from $s_\lozenge$ to $s_\ast$ in the direction $\sigma$.
\begin{lem}
\label{lem.allz} The vehicle cannot pass more than $N_s$ modes
$\mathfrak{A}$ in a row, constantly not finding the target in
$\mathfrak{C}_{\mathfrak{A}}$ and not changing the direction of the
motion along $\partial D$.
\end{lem}
\pf Suppose the contrary and that $\sigma=+$; the case $\sigma=-$ is
considered likewise. By Observation~\ref{obs.exit}, the $i$th mode
$\mathfrak{A}^i$ in the row starts when $s$ lies in an $+$exit arc
$A^i$, whereas $\zeta \geq 0$ when it ends. Hence $A^1,A^2, \ldots$
cannot repeat until $s$ completes the full run over $\partial D$.
However, they do repeat since the number of $+$arcs does not exceed
$F$ by Observation~\ref{obs.exit}, and $F \leq N_s$ by construction
from the proof of Lemma~\ref{lem.ns}. Hence the path $\mathscr{P}$
can be truncated so that the first and last modes $\mathfrak{A}$
start at positions $\bldr_1$ and $\bldr_2$, respectively, lying on a
common $+$exit arc $A$, whereas $s$ encircles the entire boundary
$\partial D$ during the move over the truncated $\mathscr{P}$. By
the definition of the $+$arc, $\fr_{\partial D}(s)$ evolves within
the fourth quadrant as $s$ runs from $s_1$ to $s_2$ within the
$+$arc and so the absolute value of its turning angle does not
exceed $\pi/2$. This and \eqref{rem.geom} (where $d_\ast:=0$) imply
that $\sphan{0}{\fr_{\partial D|_{s_1 \to s_2}}} \leq -3/2 \pi$. In
\eqref{basic.eq}, $|\sphan{\bt}{[ \bldr_i, s_i]}| < \pi/2$ and
$\angle \left[ T_{\partial D}(s_i),
T_{\mathscr{P}}(\bldr_i)\right]=0$ since the segments $[ \bldr_i,
s_i]$ and $[\bldr_i, \bt]$ are perpendicular. Overall,
\eqref{basic.eq} implies that
\begin{equation}
\label{eq.consist} \sphan{0}{\fr_{\mathscr{P}|_{\bldr_1 \to
\bldr_2}}} < - \frac{\pi}{2}.
\end{equation}
The path $\mathscr{P}|_{\bldr_1 \to \bldr_2}$ starts with $\beta=0$
and whenever $\beta =0$ is encountered, the angle $\beta$ may stay
constant during SMT but after this SMT $\beta$ becomes positive by
\eqref{eq.cross} (see Fig.~\ref{fig.cross}(b)) since the robot turns
right. The last claim holds thanks to (iii) of
Proposition~\ref{lem.cir} if $\mathfrak{B}$ is not terminated during
this SMT and \eqref{sigma.update} otherwise. Such behavior of
$\beta$ is inconsistent with \eqref{eq.consist}. The contradiction
obtained completes the proof. \epf
\par
{\it Proof of Proposition}~\ref{prop.det} is straightforward from
(v) of Proposition~\ref{lem.cir} and Lemmas~\ref{lem.find},
\ref{lem.direct}, and \ref{lem.allz}.
\subsection{Proof of (ii) in Theorem~\ref{th.main}.} \label{subsec.prth} Let $P_k$ be the probability
that the vehicle does not arrive at $\bt$ after making $kN$ switches
$\mathfrak{A} \to \mathfrak{B}$, where $N$ is taken from
Proposition~\ref{prop.det}. Given a realization of $\sigma$'s for
the first $kN$ switches, the probability of the $(k+1)$th event does
not exceed the probability $P_\ast$ that the next $N$ realizations
are not identical to those generated by the algorithm $\mathscr{A}$
for the related initial state. Here $P_\ast \leq \rho$, where $\rho
:= 1-\min\{p, 1-p\}^N$ and $p$ is the probability of picking $+$ in
\eqref{c.a}. So the law of total probability yields that $P_{k+1}
\leq \rho P_k \Rightarrow P_k \leq \rho^{k-1} P_1 \to 0$ as $k \to
\infty$. It remains to note that the probability not to achieve
$\bt$ does not exceed $P_k$ for any $k$.
\section{Proof of (ii) in Theorem~\ref{th.main} and Theorem~\ref{th.maina}}
For the definiteness, we consider the vehicle driven by the basic algorithm with the right turns. So in any SMEC the vehicle has the obstacle to the left. The proof basically follows that from the previous section and employs many facts established there. The difference is that now we do not need to introduce an auxiliary deterministic algorithm since the examined one is deterministic itself.
\par
As before, we first consider another obstacle $\mathscr{D} \not\ni \bt$ satisfying Assumption~\ref{ass.finite1}.
Let a point $\bldr$ moves in the plane according to
the following rules:
\begin{enumerate}[{\bf r.1)}]
\item If $\bldr \not\in \mathscr{D}$, $\bldr$ moves
to $\bt$ in a straight line; $\bldr(0) \not\in \mathscr{D}$;
\item If $\bldr$ hits $C:=\partial \mathscr{D}$, it turns right and
then moves in the positive direction along the boundary, counting the
angle $\beta$;
\item This motion lasts until $\beta=0$ and new SMT is
possible;
\item The point halts as soon as it arrives at the target.
\end{enumerate}
The path traced by $\bldr$ is called the {\em symbolic path (SP)}.
Any SMT according to {\bf r.1)} except
for the first one starts and ends at some points $s_\lozenge, s_\ast
\in C$, which cut out a cave $\mathfrak{C}[s_\lozenge, s_\ast]$.
\par
We start with noting that the following specification of Observation~\ref{obs.arc} now holds.
\begin{observation}
\label{obs.arc11} As $\bldr$ moves over a
$\pm$-arc of $C$, we have $\pm \dot{\varphi}
\geq 0$. Non-singular points of $\pm$-arc lie above/below
$\mathscr{D}$.
\end{observation}
Lemma~\ref{lem.countbeta} evidently remains valid, whereas Lemma~\ref{lem.beta} holds in the following specified form.
\begin{lem}
\label{lem.beta11} Whenever SP lies on $C$, we have $\beta \geq 0$.
\end{lem}
It is easy to see by inspection that Lemma~\ref{lem.posbeta} remains true as well, where in the case from Fig.~\ref{fig.frsing} the right turn at the point $s_{\#}$ is justified by not the absence of the target in the cave but the very algorithm statement.
The following claim is analog of Lemma~\ref{lem.notchange}
\begin{lem}
\label{lem.notchange11} Suppose that after SMT starting and ending at
the points $s_\lozenge$ and $s_\ast$, respectively, SP
enters the cave $\mathfrak{C}[s_\lozenge,s_\ast]$. The this cave
contains the entire path traced before SMT at hand.
\end{lem}
\pf The proof is be retracing the arguments from the proof of Lemma~\ref{lem.notchange} with the only alteration: the point cannot enter the cave through $s_{\lozenge}$ since this violates the always positive direction of motion along the boundary. \epf
\par
Now we revert to the vehicle at hand and show that Lemma~\ref{lem.notchange11} extends on the directed path $\mathscr{P}$ traced by this vehicle.
The next lemma employs the notations $\mathfrak{A}^\dagger$ and $\sigma_{\mathfrak{A}^\dagger}$ introduced at the beginning of subsection~\ref{subsec.cl}.
\setcounter{thm}{7}
\begin{lem}\hspace{-5.0pt}$^a$
For any occurrence $\mathfrak{A}^\dagger$ of mode $\mathfrak{A}$
that holds between two modes $\mathfrak{B}$, we have $\sigma_{\mathfrak{A}^\dagger} =+$.
\end{lem}
\setcounter{thm}{27}
\pf Suppose to the contrary that $\sigma_{\mathfrak{A}^\dagger} =-$. Then according to the 'only-right-turns' option of the algorithm, the vehicle enters the cave $\mathfrak{C}_{\mathfrak{A}^\dagger}$ after termination of $\mathfrak{A}^\dagger$.
We are going to show that then similar to Lemma~\ref{lem.notchange}, this cave contains the entire path passed by the vehicle until this moment and so its initial location. Due to the first relation from \eqref{far.enough}, the last claim implies that the initial location $\bldr_0$ is also contained by a cave of $N(\dtrig)$, in violation of the assumptions of Theorem~\ref{th.maina}. This contradiction will complete the proof.
\par
Thus it remains to show that $\mathfrak{C}_{\mathfrak{A}^\dagger}$ does contain the path traced so far. Suppose the contrary. Since in the mode $\mathfrak{B}$ preceding $\mathfrak{A}^\dagger$, the vehicle has the obstacle to the left, it passes to $\mathfrak{A}^\dagger$ from inside the cave. It follows that the moment after $\mathfrak{A}^\dagger$ is not the first time when the vehicle enters the cave. Let us consider the last of these 'preceding' enters and the path $\mathscr{P}$ traced by the vehicle since this moment until the commencement of $\mathfrak{A}^\dagger$. By combining Lemma~\ref{lem.posbeta} with the arguments from the proof of Lemma~\ref{lem.find}, we conclude that this path is single and $\beta>0$ at its end, which makes mode $\mathfrak{A}^\dagger$ impossible. The contradiction obtained completes the proof. \epf
\par
This lemma entails that Corollaries~\ref{lem.turn1}, \ref{lem.parr1}, and \ref{cor.lasrt} remain true in the following specified forms.
\begin{corollary}
\label{lem.turn11} \!For $\bldr_1=
\bldr_\lozenge, \bldr_2= \bldr_\ast$, \eqref{basic.eq} holds with $\sigma\!=\!+$.
\end{corollary}
\begin{corollary}
\label{lem.parr11} Let $\bldr_1, \bldr_2$ be successively passed
within a common mode $\mathfrak{B}$. If $\bldr_2$ follows IT, \eqref{basic.eq}
holds with $\sigma=+$ and $\sphan{0}{\fr_{\partial D}|_{s_1 \xrightarrow{\sigma}
s_2}}$ accounting for possible full runs over $C$.
\end{corollary}
\begin{corollary}
Suppose that points $\bldr_1$ and $\bldr_2$ are successively passed in modes
$\mathfrak{B}$ (maybe, different) and $\bldr_2$ is not
attributed to IT. Then
\eqref{basic.eq} holds with $\sigma=+$, where $\sphan{0}{\fr_{\partial D}|_{s_1
\xrightarrow{\sigma} s_2}}$ accounts for the entire motion of the
projection $s=s[\bldr], \bldr \in \mathscr{P}_{\bldr_1 \to
\bldr_2}$, including possible full runs over $\partial D$.
\end{corollary}
We also note that at the moment when a SMEC ends, $s \in S_0:= \{s
\in \partial D: -\dtrig \leq \zeta_{\partial D}(s) <0,
\lambda_{\partial D}(s) >0\}$. Since the boundary is piece-wise analytical,
this set has finitely many connected components (called exit arcs).
\par
{\bf PROOF OF THEOREM~\ref{th.maina}} This proof retraces many arguments from the proof of Lemma~\ref{lem.allz}.
Suppose the contrary and that the vehicle does not arrive at the target.
Then the projection $s$ repeatedly encircles the boundary. (This includes the imaginary moves of $s$ when the vehicle is in mode $\mathfrak{A}$.)
By retracing the arguments from the proof of (v) in Proposition~\ref{lem.cir}, we conclude that the path $\mathscr{P}$
can be truncated so that the first and last modes $\mathfrak{A}$
start at positions $\bldr_1$ and $\bldr_2$, respectively, lying on a
common exit arc $A$, whereas $s$ encircles the entire boundary
$\partial D$ during the move over the truncated $\mathscr{P}$. By
the definition of the exit arc, $\fr_{\partial D}(s)$ evolves within
the fourth quadrant as $s$ runs from $s_1$ to $s_2$ within the
$+$arc and so the absolute value of its turning angle does not
exceed $\pi/2$. This and \eqref{rem.geom} (where $d_\ast:=0$) imply
that $\sphan{0}{\fr_{\partial D|_{s_1 \to s_2}}} \leq -3/2 \pi$. In
\eqref{basic.eq}, $|\sphan{\bt}{[ \bldr_i, s_i]}| < \pi/2$ and
$\angle \left[ T_{\partial D}(s_i),
T_{\mathscr{P}}(\bldr_i)\right]=0$ since the segments $[ \bldr_i,
s_i]$ and $[\bldr_i, \bt]$ are perpendicular. Overall,
\eqref{basic.eq} implies \eqref{eq.consist}.
The path $\mathscr{P}|_{\bldr_1 \to \bldr_2}$ starts with $\beta=0$
and whenever $\beta =0$ is encountered, the angle $\beta$ may stay
constant during SMT but after this SMT $\beta$ becomes positive since the robot turns
right. The last claim holds thanks to (iii) of
Proposition~\ref{lem.cir} if $\mathfrak{B}$ is not terminated during
this SMT and the right-turn option in \eqref{c.a} otherwise. Such behavior of
$\beta$ is inconsistent with \eqref{eq.consist}. The contradiction
obtained completes the proof. \epf
\par
{\bf PROOF OF (ii) IN THEOREM~\ref{th.main}} This claim is immediate from Theorem~\ref{th.maina}. \epf


\end{document}